\documentclass[reqno]{tran-l} 

\usepackage[all]{xy}
\usepackage{rotating}
\usepackage{supertabular}
\def\theoremnumberstyle{section} % legt Stil der Satznumerierung fest
\def\proofname{Proof}
\usepackage{tom}

\newcommand{\tom}[1]{}   % Anmerkungen, die manchmal ausgeklammert werden sollen.
\newcommand{\tomtom}[1]{}  % Anmerkungen, die meistens ausgeklammert werden sollen.
\newcommand{\klammer}[1]{} %% Manchmal stoeren geoeffnete Klammern
\renewenvironment{equationlist}{
  \begin{list}{}
    {\renewcommand{\makelabel}[1]{\stepcounter{equation}\mbox{\rm(\theequation)}%
        \if##1\empty\else\mylabel{##1}{\theequation}\fi}
     \leftmargin2cm \itemindent0cm\labelwidth1.5cm \parsep6pt \topsep6pt\labelsep0.2cm}
    }
  {\end{list}
  }
\newenvironment{mylist}{
  \begin{list}{}
    {\leftmargin4ex \itemindent0ex\labelwidth2ex \parsep0.5ex \topsep0pt\labelsep1ex}
    }
  {\end{list}
  }
\newenvironment{myenumerate}{
  \begin{list}{}
    {\leftmargin5ex \itemindent0ex\labelwidth4ex \parsep0.5ex \topsep0pt\labelsep1ex}
    }
  {\end{list}
  }

\numberwithin{equation}{section} %% Reset equation counter to zero at
\allowdisplaybreaks[3]   %% Abgesetzte Formeln sollen am Seitenende
\newcommand{\asf}{1.4}                                %% Arraystrech-Faktor
\renewcommand{\arraystretch}{\asf}                    %% Arraystrech-Faktor setzen
\newcommand{\as}[1]{\renewcommand{\arraystretch}{#1}} %% Arraystrech-Faktor
\newcommand{\expl}[1]{{_{_{\mbox{\tiny #1}}}}}  % tiefgestellte Erklaerungen
\SelectTips{cm}{12}  % sets the arrow style in xymatrices to the cm-standard
\CompileMatrices     % Diagramme sollen nur nach Aenderungen neu compiliert werden

\begin{document}

   \parindent0cm
   \thispagestyle{empty}

%   \begin{titlepage}
     \title[Irreducibility]{Irreducibility of Equisingular Families of Curves}
     \author{Thomas Keilen}
   \address{Universit\"at Kaiserslautern\\
     Fachbereich Mathematik\\
     Erwin-Schr\"odinger-Stra\ss e\\
     D -- 67663 Kaiserslautern
     }
   \email{keilen@mathematik.uni-kl.de}
   \urladdr{http://www.mathematik.uni-kl.de/\textasciitilde keilen}
   \thanks{The author was partially supported by the
     DFG-Schwerpunkt ``Globale Methoden in der komplexen Geometrie''.
     The author would like to thank the referee for pointing out Example \ref{ex:referee}.}
%   \thanks{The author would like to express his thanks to Gert-Martin
%     Greuel, Christoph Lossen, and Eugenii Shustin for many helpful discussions.}


   \date{July, 2001.}

     
   \begin{abstract}
     In 1985 Joe Harris (cf.~\cite{Har85}) proved the long standing
     claim of Severi that equisingular families of plane nodal curves
     are irreducible whenever they are non-empty. For families with
     more complicated singularities this is no longer true. Given a
     divisor $D$ on a smooth projective surface $\Sigma$ it thus
     makes sense to look for conditions which ensure that the family
     $V_{|D|}^{irr}\big(\ks_1,\ldots,\ks_r\big)$ of irreducible curves
     in the linear system $|D|_l$ with precisely $r$ singular points
     of types $\ks_1,\ldots,\ks_r$ is irreducible. Considering
     different surfaces including general surfaces in
     $\PC^3$ and products of curves, we produce a sufficient condition
     of the type
     \begin{displaymath}
       \sum\limits_{i=1}^r\deg\big(X(\ks_i)\big)^2
       <
       \gamma\cdot (D- K_\Sigma)^2,       
     \end{displaymath}
     where $\gamma$ is some constant and $X(\ks_i)$ some
     zero-dimensional scheme associated to the singularity type. Our
     results carry the same asymptotics as the best known results in
     this direction in the plane case, even though the coefficient is
     worse (cf.~\cite{GLS00}). For most of the considered surfaces these
     are the only known results in that direction. 
   \end{abstract}

   \maketitle

   \tableofcontents

%% Einleitung

   \section{Introduction}

   Equisingular families of curves have been studied quite intensively since the
   last century. If we fix a linear system $|D|_l$ on a smooth
   projective surface $\Sigma$ and singularity types
   $\ks_1,\ldots,\ks_r$ we denote by
   $V^{irr}=V_{|D|}^{irr}\big(\ks_1,\ldots,\ks_r\big)$ the variety 
   of irreducible curves in $|D|_l$ with
   precisely $r$ singular points of the given types. The main
   questions are whether the equisingular family $V^{irr}$ is non-empty, smooth of the
   expected dimension, and irreducible. For results in the plane case we
   refer to \cite{GLS98a,GLS00}, and  results on the first and the
   second question on other surfaces may be found in
   \cite{GLS97,GLS98b,CC99,Fla01,Che01,KT02}. In this paper for the first time  the
   question of the irreducibility of $V^{irr}$ for a wider range of
   surfaces is studied. As already families of cuspidal curves in the plane (cf.~\cite{Zar35}) or
   nodal curves on surfaces in $\PC^3$ (cf.~\cite{CC99}) show,
   in general we cannot expect a complete answer as for families of
   plane nodal curves, saying that the family is irreducible whenever
   it is non-empty. All we may hope for are numerical conditions
   depending on invariants of the singularity types, the surface and
   the linear system, which ensure the irreducibility of $V^{irr}$. 
   
   The main condition which we get (cf.~Section \ref{sec:irred}) looks like
   \begin{equation}\label{eq:irred:1}
     \sum\limits_{i=1}^r\deg\big(X(\ks_i)\big)^2
     <
     \gamma\cdot (D- K_\Sigma)^2,     
   \end{equation}
   where $\gamma$ is some  constant. Applying the estimates
   \eqref{eq:deg:1} for
   $\deg\big(X(\ks_i)\big)$ from Subsection \ref{subsec:schemes} we could replace
   \eqref{eq:irred:1} by
   \begin{equation}\label{eq:irred:2}
     \sum\limits_{i=1}^r\tau(\ks_i)^2
     <
     \tfrac{\gamma}{9}\cdot (D- K_\Sigma)^2,     
   \end{equation}
   in the case of analytical types, and in the topological case by
   \begin{equation}\label{eq:irred:3}
     \sum\limits_{i=1}^r\big(\mu(\ks_i)+\tfrac{4}{3}\big)^2
     <
     \tfrac{4\cdot\gamma}{9}\cdot (D- K_\Sigma)^2.  
   \end{equation}

   In this section we introduce the basic concepts and notations used
   throughout the paper, and we state several important known
   facts. Section \ref{sec:irred} contains the main results and their
   proofs, omitting the technical details. These are presented in 
   Section \ref{sec:v-irr-reg} and Section
   \ref{sec:technical-lemmata}. 

%%%%%%%%%%%%%%%%%%%%%%%%%%%%%%%%%%%%%%%%%%%%%%%%%%%%%%%%%%%%%%%%%%%%%%%%%%%%%%

   \subsection{General Assumptions and Notations}\label{subsec:notations}
     Throughout this article $\Sigma$ will denote a smooth projective surface
     over $\C$. $\N$ denotes the set of non-negative integers.
        
     We will denote by $\Div(\Sigma)$ the 
     group of divisors on $\Sigma$ and by $K_\Sigma$ its canonical
     divisor. 
     If $D$ is any
     divisor on $\Sigma$, $\ko_\Sigma(D)$ shall be the corresponding invertible
     sheaf and we will sometimes write $H^\nu(X,D)$
     instead of $H^\nu\big(X,\ko_X(D)\big)$.  
     A \emph{curve} $C\subset\Sigma$ will be an effective (non-zero) divisor, that
     is a one-dimensional locally principal scheme, not necessarily
     reduced; however, an \emph{irreducible curve} shall be
     reduced by definition.
     $|D|_l$ denotes the
     system of curves linearly equivalent to $D$.%, while we use the notation
%     $|D|_a$ for the system of curves algebraically equivalent to $D$
%     (cf.~\cite{Har77} Ex.~V.1.7), that is, $|D|_a$ is the reduction of the
%     connected component of $\Hilb_\Sigma$, 
%     the Hilbert scheme of all curves on $\Sigma$, 
%     containing any curve algebraically equivalent to $D$ (cf.~\cite{Mum66} Chapter 
%     15).\tom{\footnote{Note that indeed the reduction of the Hilbert
%         scheme gives the Hilbert scheme $\Hilb_\Sigma^{red}$ of
%         curves on $\Sigma$ over reduced base spaces.}}
     We will use the notation
     $\Pic(\Sigma)$ for the \emph{Picard group} of $\Sigma$, that is
     $\Div(\Sigma)$ modulo linear equivalence (denoted by $\sim_l$), and
     $\NS(\Sigma)$ for the 
     \emph{N\'eron--Severi group}, that is $\Div(\Sigma)$ modulo algebraic
     equivalence (denoted by $\sim_a$).
     Given a reduced curve $C\subset\Sigma$ we will write $g(C)$ for
     its \emph{geometric genus}.
   
     Given any closed subscheme
     $X$ of a scheme $Y$, we denote by 
     $\kj_X=\kj_{X/Y}$ the  
     \emph{ideal sheaf} of $X$ in $\ko_Y$. If $X$ is zero-dimensional we denote
     by $\# X$  the number of points in its
     \emph{support}  $\supp(X)$ and
     by
     $\deg(X)=\sum_{z\in Y}\dim_\C(\ko_{Y,z}/\kj_{X/Y,z})$ 
     its \emph{degree}. 

     If $X\subset \Sigma$ is a zero-dimensional
     scheme on $\Sigma$ and $D\in\Div(\Sigma)$, we denote by
     $\big|\kj_{X/\Sigma}(D)\big|_l$ the linear system of curves $C$ in
     $|D|_l$ with $X\subset C$.
     
     If $L\subset\Sigma$ is any reduced curve and $X\subset \Sigma$ a zero-dimensional scheme,
     we define the \emph{residue scheme}
     $X:L\subset\Sigma$ \emph{of $X$} by the ideal sheaf  
     $\kj_{X:L/\Sigma}=\kj_{X/\Sigma}:\kj_{L/\Sigma}$ with stalks
     \begin{displaymath}
       \kj_{X:L/\Sigma,z}= \kj_{X/\Sigma,z}:\kj_{L/\Sigma,z},
     \end{displaymath}
     where ``$:$'' denotes the ideal quotient. This  leads to the
     definition of the \emph{trace scheme} $X\cap L\subset L$ \emph{of
       $X$} via the 
     ideal sheaf $\kj_{X\cap L/L}$ given by the exact sequence
     \begin{displaymath}
       \xymatrix@C0.6cm{
         0\ar[r] & {\kj_{X:L/\Sigma}(-L)}\ar[r]^(0.6){\cdot L} & {\kj_{X/\Sigma}}\ar[r]
         &{\kj_{X\cap L/L}}\ar[r] &0.
         }
     \end{displaymath}

   \subsection{Singularity Types}\label{subsec:types}
     The germ $(C,z)\subset(\Sigma,z)$ of a reduced curve $C\subset\Sigma$ at
     a point $z\in\Sigma$ is called a \emph{plane curve singularity},
     and two plane curve singularities $(C,z)$ and $\big(C',z'\big)$ are said to be
     \emph{topologically} (respectively \emph{analytically equivalent})
     if there is homeomorphism (respectively an analytical isomorphism) 
     $\Phi:(\Sigma,z)\rightarrow(\Sigma,z')$ such that
     $\Phi(C)=C'$. We call an equivalence class with respect to these
     equivalence relations a \emph{topological} (respectively
     \emph{analytical singularity type}). 
     The following are known to be invariants of the topological type $\ks$
     of the plane curve singularity $(C,z)$:
     $r(\ks)=r(C,z)$, the number of branches of
     $(C,z)$;
     $\tau^{es}(\ks)=\tau^{es}(C,z)$, the codimension of the $\mu$-constant
     stratum in the semiuniversal deformation of
     $(C,z)$; 
     $\delta(\ks)=\delta(C,z)=\dim_\C\big(\nu_*\ko_{\widetilde{C},z}/\ko_{C,z}\big)$, the 
     \emph{delta invariant} of $\ks$, where
     $\nu:\big(\widetilde{C},z\big)\rightarrow(C,z)$ is a
     normalisation of $(C,z)$; and
     $\mu(\ks)=\mu(C,z)=\dim_\C\ko_{\Sigma,z}\big/\big(\frac{\partial
       f}{\partial x},\frac{\partial f}{\partial y}\big)$, the
     \emph{Milnor number} of $\ks$, where
     $f\in\ko_{\Sigma,z}$ denotes a local equation of $(C,z)$ with respect to the
     local coordinates $x$ and $y$.
     For
     the analytical type $\ks$ of $(C,z)$ we have as additional invariant the
     \emph{Tjurina number} of $\ks$ defined as
     $\tau(\ks)=\tau(C,z)=\dim_\C\ko_{\Sigma,z}\big/\big(f,\frac{\partial
       f}{\partial x},\frac{\partial f}{\partial y}\big)$.
     We recall the relation $2\delta(\ks)=\mu(\ks)+r(\ks)-1$
     (cf.~\cite{Mil68}   Chapter 10). 
     Furthermore, since the
     $\delta$-constant stratum of the semiuniversal deformation of
     $(C,z)$ contains the $\mu$-constant stratum and
     since its codimension is just $\delta(\ks)$, we have
     $\delta(\ks)\leq\tau^{es}(\ks)$ (see also \cite{DH88}); and hence
     \begin{equation}\label{eq:types}
       \mu(\ks)
       \leq
       2\delta(\ks)\leq 2\tau^{es}(\ks).
     \end{equation}

   \subsection{Singularity Schemes}\label{subsec:schemes}          
     For a reduced curve $C\subset\Sigma$ we recall the definition of
     the zero-dimensional schemes $X^{es}(C)\subseteq X^s(C)$ and
     $X^{ea}(C)\subseteq X^a(C)$ from \cite{GLS00}. They are defined
     by the ideal sheaves $\kj_{X^{es}(C)/\Sigma}$,
     $\kj_{X^s(C)/\Sigma}$, $\kj_{X^{ea}(C)/\Sigma}$, and
     $\kj_{X^a(C)/\Sigma}$ respectively, given by the following stalks
     \begin{mylist}
     \item[$\bullet$] $\kj_{X^{es}(C)/\Sigma,z}=I^{es}(C,z)=
       \big\{g\in\ko_{\Sigma,z}\;\big|\;f+\varepsilon g \mbox{ is equisingular
         over } \C[\varepsilon]/(\varepsilon^2)\big\}$,
       where $f\in\ko_{\Sigma,z}$ is a local equation of $C$ at $z$.
       $I^{es}(C,z)$ is called the \emph{equisingularity ideal} of
       $(C,z)$.
     \item[$\bullet$] $\kj_{X^{s}(C)/\Sigma,z}=%I^{s}(C,z)=
       \Big\{g\in\ko_{\Sigma,z}\;\Big|\; g \mbox{ goes through the cluster }
       \Cl\big(C,T^*(C,z)\big)\Big\}$,
       where $T^*(C,z)$ denotes the essential subtree of the complete
       embedded resolution tree of $(C,z)$. 
%       $I^{s}(C,z)$ is called the \emph{singularity ideal} of $(C,z)$.
     \item[$\bullet$] $\kj_{X^{ea}(C)/\Sigma,z}=I^{ea}(C,z)=% I^{ea}(f)=
       \big( f,\tfrac{\partial f}{\partial x},\tfrac{\partial
         f}{\partial y}\big)\subseteq\ko_{\Sigma,z}$, 
       where $x,y$ denote local coordinates of $\Sigma$ at $z$ and
       $f\in\ko_{\Sigma,z}$ is a local equation of $C$.
       $I^{ea}(C,z)$ is called the \emph{Tjurina ideal} of $(C,z)$.
     \item[$\bullet$] $\kj_{X^a(C)/\Sigma,z}=I^a(C,z)
       \subseteq\ko_{\Sigma,z}$, 
       where we refer  for the somewhat lengthy definition of
       $I^a(C,z)$ to \cite{GLS00} Section 1.3.% Remark~\ref{rem:Ia}.        
%       $I^{a}(C,z)$ is called the \emph{analytical singularity
%         ideal} of the  
%       singularity $(C,z)$.
     \end{mylist}
     We call $X^{es}(C)$ the \emph{equisingularity scheme} of $C$ and
     $X^s(C)$ its \emph{singularity scheme}. Analogously we call
     $X^{ea}(C)$ the \emph{equianalytical singularity scheme} of
     $C$ and $X^a(C)$ its \emph{analytical singularity scheme}.
     \begin{center}
       \framebox[11cm]{
         \begin{minipage}{10cm}
         \medskip
           Throughout this article we will frequently treat
           topological and analytical singularities at the same time.
           Whenever we do so, we will write $X^*(C)$ for $X^{es}(C)$
           respectively for $X^{ea}(C)$ and similarly $X(C)$ for
           $X^s(C)$ respectively for $X^a(C)$.
         \medskip
         \end{minipage}
         }
     \end{center}

     In \cite{Los98}, Propositions 2.19 and 2.20 and in Remarks 2.40 (see
     also \cite{GLS00}) and
     2.41, it is shown that, fixing a point $z\in\Sigma$ and a
     topological (respectively analytical) type $\ks$, the singularity schemes
     (respectively analytical) singularity schemes having the same topological (respectively
     analytical) type are 
     parametrised by an irreducible Hilbert scheme, which we are going
     to denote by $\Hilb_z(\ks)$. This then leads to an irreducible
     family  
     \begin{equation}\label{eq:dim-hilb}
       \Hilb(\ks)=\coprod_{z\in\Sigma}\Hilb_z(\ks).
     \end{equation}
     In particular, equisingular (respectively equianalytical)
     singularities have singularity schemes (respectively analytical singularity schemes)
     of the same degree (see also \cite{GLS98a} or \cite{Los98} Lemma 2.8). The same is
     of course true, regarding the equisingularity scheme (respectively the
     equianalytical singularity scheme). 
     If $C\subset\Sigma$ is a reduced curve such that $z$ is a
     singular point of topological (respectively analytical) type $\ks$,
     we may therefore define $\deg\big(X(\ks)\big)=\deg\big(X(C),z\big)$
     and $\deg\big(X^{*}(\ks)\big)=\deg\big(X^{*}(C),z\big)$.
     We note that, with this notation,
     $\dim\Hilb_z(\ks)=\deg\big(X(\ks)\big)-\deg\big(X^*(\ks)\big)-2$
     for any $z\in\Sigma$, and thus 
     \begin{displaymath}
       \dim\Hilb(\ks)=\deg\big(X(\ks)\big)-\deg\big(X^*(\ks)\big).
     \end{displaymath}
     In the applications it is convenient to replace
     the degree of an (analytical) singularity scheme by an upper
     bound in known invariants of the singularities. From \cite{Los98}
     p.~28, p.~103, and Lemma 2.44 it follows for a topological
     (respectively analytical) singularity type $\ks$ one has
     \begin{equation}\label{eq:deg:1}
       \deg\big(X^a(\ks)\big)\leq 3\tau(\ks) \;\;\;\text{ and } \;\;\;
       \deg\big(X^s(\ks)\big)\leq \tfrac{3}{2}\mu(\ks)+2.
     \end{equation} 

   \subsection{Equisingular Families}\label{subsec:families}
     Given a divisor $D\in\Div(\Sigma)$ and topological or analytical singularity types
     $\ks_1,\ldots, \ks_r$, we denote
     by $V=V_{|D|}(\ks_1,\ldots,\ks_r)$ the locally closed subspace of
     $|D|_l$ of reduced curves in the linear system $|D|_l$ having
     precisely $r$ singular points of types $\ks_1,\ldots,\ks_r$.      
     By\footnote{$V^{reg}$ should not be confused with $\big\{C\in
       V\;\big|\;h^1\big(\Sigma,\kj_{X^*(C)/\Sigma}(D)\big)=0\big\}$,
       which is the part of $V$, where $V$ is smooth of the expected
       dimension. Curves in the latter subscheme are often called
       \emph{regular} (c.\ f.\ \cite{CC99}). See also Example
       \ref{ex:referee}.} 
     $V^{reg}=V_{|D|}^{reg}(\ks_1,\ldots,\ks_r)$ we denote the
     open (cf.~Proof of Theorem \ref{thm:v-reg}) subset  
     \begin{displaymath}
       V^{reg}=\big\{C\in V\;\big|\;h^1\big(\Sigma,\kj_{X(C)/\Sigma}(D)\big)=0\big\}
       \subseteq V.
     \end{displaymath}
     Similarly, we use the notation
     $V^{irr}=V_{|D|}^{irr}(\ks_1,\ldots,\ks_r)$ to denote the open
     subset of irreducible curves in the space $V$, and we set
     $V^{irr,reg}=V_{|D|}^{irr,reg}(\ks_1,\ldots,\ks_r)=V^{irr}\cap
     V^{reg}$, which is  open in $V^{reg}$ and in
     $V$.    
     If a type $\ks$ occurs $k>1$ times, we rather write $k\ks$
     than $\ks,\stackrel{k}{\ldots},\ks$.
     We call these families of curves \emph{equisingular families of
       curves}.

     We say that $V$ is \emph{T-smooth} at $C\in V$ if the germ
     $(V,C)$ is smooth of the (expected) dimension
     $\dim|D|_l-\deg\big(X^*(C)\big)$.

     By \cite{Los98}
     Proposition 2.1 (see also \cite{GK89}, \cite{GL96}, \cite{GLS00})
     T-smoothness of $V$ at $C$ follows from the vanishing of

     $H^1\big(\Sigma,\kj_{X^*(C)/\Sigma}(C)\big)$. This is
     due to the fact that the tangent space of $V$ at $C$ may be
     identified with
     $H^0\big(\Sigma,\kj_{X^*(C)/\Sigma}(C)\big)/H^0(\Sigma,\ko_\Sigma)$.

   \subsection{Fibrations}\label{subsec:psi}
     Let $D\in\Div(\Sigma)$ be a divisor, $\ks_1,\ldots,\ks_r$
     distinct topological or analytical singularity types, and $k_1,\ldots,
     k_r\in\N\setminus\{0\}$. 
     We denote by $\widetilde{B}$ the
     irreducible parameter space 
     \begin{displaymath}
       \widetilde{B}=\widetilde{B}(k_1\ks_1,\ldots,k_r\ks_r)
       =\prod_{i=1}^r\Sym^{k_i}\big(\Hilb(\ks_i)\big),
     \end{displaymath}
     and by $B=B(k_1\ks_1,\ldots,k_r\ks_r)$  
     the non-empty open, irreducible and dense subspace
     \begin{multline*}
         B=\Big\{\big([X_{1,1},\ldots,X_{1,k_1}],\ldots,[X_{r,1},\ldots,X_{r,k_r}]\big)
           \in\widetilde{B}\;\Big|\;\supp(X_{i,j})\cap\supp(X_{s,t})=\emptyset\\
         \;\;\forall\;1\leq i,s\leq r, 1\leq j\leq k_i, 1\leq t\leq
         k_s\Big\}.
     \end{multline*}
     Note that $\dim(B)$ does not depend on $\Sigma$; more
     precisely, with the notation of Subsection \ref{subsec:schemes} 
     we have
     \begin{displaymath}
       \dim(B)=\sum_{i=1}^r k_i\cdot \Big(\deg\big(X(\ks_i)\big)-\deg\big(X^*(\ks_i)\big)\Big).
     \end{displaymath}

     Let us set $n=\sum_{i=1}^r k_i\deg\big(X(\ks_i)\big)$. We then
     define an injective morphism  
     \begin{displaymath}
       \xymatrix@R=0.2cm{
         \psi=\psi(k_1\ks_1,\ldots,k_r\ks_r):
         B(k_1\ks_1,\ldots,k_r\ks_r)\ar[r]&
         \Hilb_\Sigma^n\\
         \;\big([X_{1,1},\ldots,X_{1,k_1}],\ldots,[X_{r,1},\ldots,X_{r,k_r}]\big)
         \ar@{|->}[r]
         &\;\bigcup_{i=1}^r\bigcup_{j=1}^{k_i}X_{i,j},
         }
     \end{displaymath}
     where $\Hilb_\Sigma^n$ denotes the smooth connected Hilbert scheme of
     zero-dimensional schemes of degree $n$ on $\Sigma$
     (cf.~\cite{Los98} Section 1.3.1). 
     
     We denote  by
     $\Psi=\Psi_D(k_1\ks_1,\ldots,k_r\ks_r)$ the fibration of
     $V_{|D|}(k_1\ks_1,\ldots,k_r\ks_r)$ induced by
     $B(k_1\ks_1,\ldots,k_r\ks_r)$; in other words
     the morphism $\Psi$ is given by  
     \begin{displaymath}
       \xymatrix@R=0.2cm{
         \Psi:
         V_{|D|}(k_1\ks_1,\ldots,k_r\ks_r)\ar[r]&
         B(k_1\ks_1,\ldots,k_r\ks_r)\\
         C\ar@{|->}[r]&
         \;\big([X_{1,1},\ldots,X_{1,k_1}],\ldots,[X_{r,1},\ldots,X_{r,k_r}]\big)
         }
     \end{displaymath}
     where $\Sing(C)=\{z_{i,j}\;|\;i=1,\ldots,r,j=1,\ldots,k_i\}$,
     $X_{i,j}=X(C,z_{i,j})$ 
     and $(C,z_{i,j})\cong\ks_i$ for all $i=1,\ldots,r,
     j=1,\ldots,k_i$. 

     With notation of Subsection \ref{subsec:families} 
     note that for $C\in V$ the fibre $\Psi^{-1}\big(\Psi(C)\big)$ is
     the open dense subset of the linear system
     $\big|\kj_{X(C)/\Sigma}(D)\big|_l$ consisting of
     the curves $C'$ with $X\big(C'\big)=X(C)$. In particular, the
     fibres of $\Psi$ restricted to $V^{reg}$ are irreducible, and since
     for $C\in V^{reg}$ the cohomology group
     $H^1\big(\Sigma,\kj_{X(C)/\Sigma}(D)\big)$ vanishes, they
     are equidimensional of dimension  
     \begin{displaymath}       
       h^0\big(\Sigma,\kj_{X(C)/\Sigma}(D)\big)-1
%       =h^0\big(\Sigma,\ko_\Sigma(D)\big)-\deg\big(X(C)\big)-1
       =h^0\big(\Sigma,\ko_\Sigma(D)\big)-\sum_{i=1}^r k_i\deg\big(X(\ks_i)\big)-1.
     \end{displaymath}     

     \tom{
       \begin{proof}[Proof of Equation \eqref{eq:dim-hilb}]
         Since $\dim\Hilb(\ks)$ is independent of $\Sigma$, we may as well
         suppose that $\Sigma=\PC^2$, and we let $H$ be a line in $\PC^2$.
         
         For a reduced curve $C\subset\PC^2$ we set $X(C)=X^s(C)$ and
         $X^*(C)=X^{es}(C)$ (respectively $X(C)=X^a(C)$ and $X^*(C)=X^{ea}(C)$).
         
         Since $\deg\big(X\big)=\deg\big(X(\ks)\big)$ is independent of $X\in\Hilb(\ks)$,
         there is an integer $m>0$ such that for $d> m$
         \begin{displaymath}
           h^1\big(\PC^2,\kj_{X/\PC^2}(d)\big)=0
         \end{displaymath}
         for any $X\in\Hilb(\ks)$.
         
         Let $k>0$ be the determinacy bound
         of $\ks$, that is, any representative
         $f\in\ko_{\PC^2,z}=\C\{x,y\}$ of $\ks$ depends only on the $k$-jet of
         $f$. Hence, for $d> k$ the morphism
         \begin{displaymath}
           \Psi=\Psi_{dH}(\ks):V_{|dH|}(\ks)\rightarrow B(\ks)=\Hilb(\ks)
         \end{displaymath}
         is surjective. 
         
         Let us now fix some $d>\max\{k,m\}$.
         For each $C\in V_{|dH|}(\ks)$ the fibre
         $\Psi^{-1}\big(\Psi(C)\big)$ is the open dense subset of
         $\big|\kj_{X(C)/\PC^2}(d)\big|$, consisting
         of curves $C'\in V_{|dH|}(\ks)$ with $X\big(C'\big)=X(C)$. 
         From the long exact cohomology sequence of 
         \begin{displaymath}
           0\rightarrow \kj_{X(C)/\Sigma}\rightarrow
           \ko_{\PC^2}(d)\rightarrow \ko_{X(C)}\rightarrow 0
         \end{displaymath}
         it follows
         \begin{displaymath}
           h^0\big(\PC^2,\kj_{X(C)/\PC^2}(d)\big)%=h^0\big(\PC^2,\ko_{\PC^2}(d)\big)-\deg\big(X(C)\big)
           =h^0\big(\PC^2,\ko_{\PC^2}(d)\big)-\deg\big(X(\ks)\big). 
         \end{displaymath}
         In particular the fibres all have the same dimension
         \begin{displaymath}
           \dim\Psi^{-1}\big(\Psi(C)\big)=h^0\big(\PC^2,\ko_{\PC^2}(d)\big)-\deg\big(X(\ks)\big)-1.
         \end{displaymath}
         We therefore get
         \begin{multline*}
           \dim\Hilb(\ks)=\dim\Psi\big(V_{|dH|}(\ks)\big)
           =\dim\big(V_{|dH|}(\ks)\big)-\dim\Psi^{-1}\big(\Psi(C)\big)\\[0.2cm]
           =\dim\big(V_{|dH|}(\ks)\big)-h^0\big(\PC^2,\ko_{\PC^2}(d)\big)+\deg\big(X(\ks)\big)+1.
         \end{multline*}
         Moreover, by Subsection \ref{subsect:schemes} we know that also
         $h^1\big(\PC^2,\kj_{X^*(C)/\PC^2}(d)\big)=0$ for any $C\in
         V_{|dH|}(\ks)$, and thus in view of 
         Subsection \ref{subsec:families} $V_{|dH|}(\ks)$ is T-smooth, that is 
         \begin{displaymath}
           \dim\big(V_{|dH|}(\ks)\big)=h^0\big(\PC^2,\ko_{\PC^2}(d)\big)-\deg\big(X^*(\ks)\big)-1,
         \end{displaymath}
         which finishes the claim.
       \end{proof}
       }

%%%%%%%%%%%%%%%%%%%%%%%%%%%%%%%%%%%%%%%%%%%%%%%%%%%%%%%%%%%%%%%%%%%%%%%%%%%%%%%%%

   \section{The Main Results}\label{sec:irred}

   In this section we give sufficient conditions for the
   irreducibility of equisingular families of curves on certain
   surfaces with Picard number one -- including the projective
   plane, general surfaces in $\PC^3$ and general K3-surfaces --, on
   products of curves, and on a subclass of geometrically ruled
   surfaces.

   \subsection{Surfaces with Picard Number One}
   
   \begin{theorem}\label{thm:irred-p3}
     Let $\Sigma$ be a surface such that
     \begin{myenumerate}
     \item[\rm(i)] $\NS(\Sigma)=L\cdot\Z$ with $L$ ample, and
     \item[\rm(ii)] $h^1(\Sigma,C)=0$, whenever $C$ is effective.
     \end{myenumerate}

     Let $D\in\Div(\Sigma)$, let $\ks_1,\ldots,\ks_r$ be pairwise distinct
     topological or analytical singularity types and let
     $k_1,\ldots,k_r\in\N\setminus\{0\}$. 

     Suppose that 
     \begin{equationlist}
     \item[eq:irred-p3:0]  \hspace*{-0.5cm}$D-K_\Sigma$ is big and nef,
     \item[eq:irred-p3:0+]  \hspace*{-0.5cm}$D+K_\Sigma$ is nef,
     \item[eq:irred-p3:1]  \hspace*{-0.5cm}$\sum\limits_{i=1}^rk_i\deg\big(X(\ks_i)\big)<\beta\cdot
       (D-K_\Sigma)^2$\;\;
       for some $0<\beta\leq
       \tfrac{1}{4}$, and 
     \item[eq:irred-p3:2] \hspace*{-0.5cm}$\sum\limits_{i=1}^rk_i\deg\big(X(\ks_i)\big)^2
       <
       \gamma\cdot (D-K_\Sigma)^2$,
       where
       $\gamma=\tfrac{\big(1+\sqrt{1-4\beta}\big)^2\cdot
         L^2}{4\cdot\chi(\ko_\Sigma)+\max\{0,2\cdot K_\Sigma.L\}+6\cdot
         L^2}$.
     \end{equationlist}
     
     Then either
     $V_{|D|}^{irr}(k_1\ks_1,\ldots,k_r\ks_r)$ is empty
     or it is irreducible of the expected dimension.
     \hfill $\Box$
   \end{theorem}

   \begin{remark}\label{rem:irred-p3}
     If we set
     \begin{displaymath}
       \gamma=\frac{36\alpha}{(3\alpha+4)^2}\;\;\;\text{ with }\;\;\; 
       \alpha=\frac{4\cdot\chi(\ko_\Sigma)+\max\{0,2\cdot K_\Sigma.L\}+6\cdot
       L^2}{L^2},
     \end{displaymath}
     then a simple calculation shows that
     \eqref{eq:irred-p3:1} becomes redundant.
     For this we have to take into account that
     $\deg\big(X(\ks)\big)\geq 3$ for any singularity type $\ks$. 
     The claim then follows with
     $\beta=\frac{1}{3}\cdot\gamma\leq \frac{1}{4}$.
     \hfill $\Box$
   \end{remark}

   We now apply the result in several special cases.

   \begin{corollary}
     Let $d\geq 3$, $L\subset\PC^2$ be a line, and 
     $\ks_1,\ldots,\ks_r$ be 
     topological or analytical singularity types.

     Suppose that 
     \begin{displaymath}
       \sum\limits_{i=1}^r \deg\big(X(\ks_i)\big)^2
       <
       \tfrac{90}{289}\cdot (d+3)^2.
     \end{displaymath}

     Then either
     $V_{|dL|}^{irr}(\ks_1,\ldots,\ks_r)$ is empty
     or it is irreducible and T-smooth.
     \hfill $\Box$
   \end{corollary}

   Many authors were concerned with the question in the case of plane
   curves with
   nodes and cusps or with nodes and one more complicated singularity
   or simply with ordinary multiple points 
   -- cf.~e.~g.~
   \cite{Sev21,AC83,Har85,Kan89a,Kan89b,Ran89,Shu91a,Shu91b,Bar93a,Shu94,Shu96a,Shu96b,Wal96,GLS98b,GLS98d,Los98,Bru99,GLS00}. 
   Using particularly designed techniques for these cases they get of
   course better results than we may expect to. 

   The best general
   results in this case can be found in \cite{GLS00} (see also \cite{Los98} Corollary
   6.1). Given a plane curve of degree $d$, omitting nodes and cusps, they get 
   \begin{displaymath}
     \sum_{i=1}^r\big(\tau^*(\ks_i)+2\big)^2\leq \tfrac{9}{10}\cdot d^2
   \end{displaymath}
   as the main irreducibility condition, where
   $\tau^*(\ks_i)=\tau(\ks_i)$ in the analytical case (respectively
   $\tau^*(\ks_i)=\tau^{es}(\ks_i)$ in the topological case). By Subsection
   \ref{subsec:types} we know that $\mu(\ks_i)\leq
   2\cdot\tau^{es}(\ks_i)$.
   Thus, in view of \eqref{eq:irred:2},
   \eqref{eq:irred:3}, \eqref{eq:types}  and of
   Theorem \ref{thm:irred-p3} we get the sufficient condition
%   \begin{displaymath}
%     \sum_{i=1}^r \tau^*(\ks_i)^2 < \tfrac{10}{289}\cdot (d+3)^2
%   \end{displaymath}
%   in the analytical case, and in the topological case
   \begin{displaymath}
     \sum_{i=1}^r\big(\tau^*(\ks_i)+\tfrac{2}{3}\big)^2< \tfrac{10}{289}\cdot (d+3)^2,
   \end{displaymath}
   which has the same asymptotics. However, the coefficients differ by a factor of about $26$.

   A smooth complete intersection surface with Picard
   number one satisfies the assumptions of Theorem
   \ref{thm:irred-p3}. Thus by the Theorem of Noether the result
   applies in particular to general surfaces in $\PC^3$.

%   \begin{corollary}
%     Let $\PC^2\not\cong\Sigma\subset\PC^N$ be a smooth complete
%     intersection of type $(d_1,\ldots,d_{N-2})$, 
%     let $H\subset\Sigma$ be a hyperplane
%     section and suppose that the Picard number of
%     $\Sigma$ is one.

%     Let $d>\kappa=\sum_{i=1}^{N-2}d_i-N-1\geq 0$,  $n=H^2=d_1\cdots d_{N-2}$ and let $\ks_1,\ldots,\ks_r$ be 
%     topological or analytical singularity types.

%     Suppose that 
%     \begin{equation}
%       \sum\limits_{i=1}^r \deg\big(X(\ks_i)\big)^2
%       <
%     \tfrac{18\cdot\big((\kappa+3)\cdot n+2\chi(\ko_\Sigma)\big)\cdot
%       n^2}{\big((3\kappa+11)\cdot n+6\chi(\ko_\Sigma)\big)^2}
%       \cdot (d-\kappa)^2.
%     \end{equation}

%     Then either
%     $V_{|dH|}^{irr}(\ks_1,\ldots,\ks_r)$ is empty
%     or it is irreducible and T-smooth.
%   \end{corollary}
%   \begin{proof}
%     It remains
%     to show $\kappa=-N-1+\sum_{i=1}^{N-2}d_i\geq 0$, then in
%     particular $dH-K_\Sigma=(d-\kappa)\cdot H$ is ample and
%     $dH+K_\Sigma=(d+\kappa)\cdot H$ is nef.
     
%     Since $\kappa\geq 2N-4-N-1=N-5$ anyway, the
%     critical situations are $N=3$ with $d_1\leq 3$, and $N=4$ with
%     $d_1=d_2=2$. In these cases the surface $\Sigma$ is either
%     $\PC^2$ or  rational with a Picard number larger than one
%     (see p.~\pageref{page:p3} and \cite{Har77} Ex.~V.4.13).
%     This finishes the claim.
%   \end{proof}

   \begin{corollary}
     Let $\Sigma\subset\PC^3$ be a smooth hypersurface of degree
     $n\geq 4$,
     let $H\subset\Sigma$ be a hyperplane
     section, and suppose that the Picard number  of
     $\Sigma$ is one.
     Let $d>n-4$ and let $\ks_1,\ldots,\ks_r$ be 
     topological or analytical singularity types.

     Suppose that 
     \begin{displaymath}
       \sum\limits_{i=1}^r \deg\big(X(\ks_i)\big)^2
       <
       \tfrac{6\cdot\big(n^3-3n^2+8n-6\big)\cdot n^2}{\big(n^3-3n^2+10n-6\big)^2}
       \cdot (d+4-n)^2,
     \end{displaymath}

     Then either
     $V_{|dH|}^{irr}(\ks_1,\ldots,\ks_r)$ is empty
     or irreducible of the expected dimension.
     \hfill $\Box$
   \end{corollary}

   We would like to thank the referee for pointing out the following
   example of reducible families $V_{|H|}^{irr}(3 A_1)$ of nodal
   curves on surfaces in $\PC^3$.

   \begin{example}\label{ex:referee}
     If $\Sigma\subset\PC^3$ is a general surface of degree $n\geq 4$,
     then 
     there is a finite number $N>1$ of $3$-tangent planes to $\Sigma$.
     However, every $3$-tangent plane cuts out an
     irreducible $3$-nodal curve on $\Sigma$, and since the Picard
     group is generated by a hyperplane section $H$, every $3$-nodal
     curve is of this form. Therefore, $V_{|H|}^{irr}(3 A_1)$ consists
     of $N$ distinct points. It is thus reducible, but smooth of the expected
     dimension 
     \begin{displaymath}
       \dim\big(V_{|H|}^{irr}(3 A_1)\big)=
       \dim|H|_l-3=0.
     \end{displaymath}
     Note that in this situation for $C\in V_{|H|}^{irr}(3 A_1)$ and
     $z\in\Sing(C)$ we have $\kj_{X(C)/\Sigma,z}=\m_{\Sigma,z}^2$ and
     thus
     \begin{displaymath}
       h^1\big(\Sigma,\kj_{X(C)/\Sigma}(H)\big)=6>0.
     \end{displaymath}
     Therefore, $V^{irr,reg}_{|H|}(3 A_1)=\emptyset$. The parameter
     space $B$ is just $\Sym^3(\Sigma)$. 
   \end{example}

%   Calculating the invariants in \eqref{eq:irred-p3-p3:1} for the
%   examples of reducible families of irreducible nodal curves on 
%   surfaces in $\PC^3$ given in the introduction on page \pageref{page:ex-p3}
%   we end up with
%   \begin{multline*}
%     \tfrac{6\cdot\big(n^3-3n^2+8n-6\big)\cdot n^2}{\big(n^3-3n^2+10n-6\big)^2}
%     \cdot (d+4-n)^2\\
%     <
%     \tfrac{6}{n-3}\cdot(d+4-n)^2
%     <
%     \tfrac{6}{n-3}\cdot\big(d^2+(4-n)\cdot d+2)\\
%     <
%     \tfrac{9n}{2}\cdot\big(d^2+(4-n)\cdot d+2)
%     =
%     9r=\sum_{i=1}^r\big(\deg X(A_1)\big)^2,
%   \end{multline*}
%   that is, our result fits with these families.

   A general K3-surface has Picard number one and in this situation,
   by the Kodaira Vanishing Theorem
   $\Sigma$ also satisfies the assumption (ii) in Theorem \ref{thm:irred-p3}.

   \begin{corollary}
     Let $\Sigma$ be a smooth K3-surface with $\NS(\Sigma)=L\cdot\Z$
     with $L$ ample and set $n=L^2$.
     Let $d>0$, $D\sim_a dL$ and let $\ks_1,\ldots,\ks_r$ be 
     topological or analytical singularity types.

     Suppose that 
     \begin{displaymath}
       \sum\limits_{i=1}^r \deg\big(X(\ks_i)\big)^2
       <
       \tfrac{54n^2+72n}{(11n+12)^2}\cdot d^2\cdot n.
     \end{displaymath}

     Then either
     $V_{|D|}^{irr}(\ks_1,\ldots,\ks_r)$ is empty
     or irreducible of the expected dimension.
     \hfill $\Box$
   \end{corollary}

%%%%%%%%%%%%%%%%%%%%%%%%%%%%%%%%%%%%%%%%%%%%%%%%%%%%%%%%%%%%%%%%%%%%%%%%%%%%%%%%%%%%%%%%%%%%

   \subsection{Products of Curves}\label{subsec:irreducibility:product-surfaces}

   If $\Sigma=C_1\times C_2$ is the product of two smooth projective
   curves, then for a general choice of $C_1$ and $C_2$ the
   N\'eron--Severi group will be generated by two fibres of the
   canonical projections, by abuse of notation also denoted by $C_1$
   and $C_2$. If both curves are elliptic, then ``general'' just means
   that the two curves are non-isogenous.

   \begin{theorem}\label{thm:irred-products-of-curves}
     Let $C_1$ and $C_2$ be two smooth projective curves of genera $g_1$
     and $g_2$  respectively with $g_1\geq g_2\geq 0$, such that for $\Sigma=C_1\times C_2$
     the N\'eron--Severi group is $\NS(\Sigma)=C_1\Z\oplus C_2\Z$.

     Let $D\in\Div(\Sigma)$ such that $D\sim_a aC_1+bC_2$ with
     $a>\max\{2g_2-2,2-2g_2\}$ and $b>\max\{2g_1-2,2-2g_1\}$, let
     $\ks_1,\ldots,\ks_r$ be pairwise distinct
     topological or analytical singularity types and
     $k_1,\ldots,k_r\in\N\setminus\{0\}$. 

     Suppose that 
     \begin{equation}
       \label{eq:irred-products-of-curves:1}
       \sum\limits_{i=1}^rk_i\deg\big(X(\ks_i)\big)^2
       \;<\;\gamma\cdot(D-K_\Sigma)^2, 
     \end{equation} 
     where $\gamma$ may be taken from the following table with $\alpha=\frac{a-2g_2+2}{b-2g_1+2}>0$.
     \begin{center}
       \tablefirsthead{\hline
         \multicolumn{1}{|c}{$g_1$}&\multicolumn{1}{|c|}{$g_2$}&$\gamma$\\\hline\hline}
       \tablehead{\hline
         \multicolumn{1}{|c}{$g_1$}&\multicolumn{1}{|c|}{$g_2$}&$\gamma$\\\hline\hline}
       \tabletail{\hline}
       \tablelasttail{\\\hline}
       \begin{supertabular}{|r|r|c|}         
         $0$&$0$&$\frac{1}{24} $\\
         $1$&$0$&$\frac{1}{\max\{32,2\alpha\}}$\\
         $\geq 2$&$0$&
         $\frac{1}{\max\{24+16g_1,4g_1\alpha\}} $\\
         $1$&$1$&
         $\frac{1}{\max\big\{32,2\alpha,\tfrac{2}{\alpha}\big\}} $\\
         $\geq 2$&$\geq 1$&
         $\frac{1}{\max\left\{24+16g_1+16g_2,4g_1\alpha,\tfrac{4g_2}{\alpha}\right\}}$
       \end{supertabular}
     \end{center}

     Then either
     $V_{|D|}^{irr}(k_1\ks_1,\ldots,k_r\ks_r)$ is empty
     or it is irreducible of the expected dimension.
     \hfill $\Box$
   \end{theorem}

   Only in the case $\Sigma\cong\PC^1\times\PC^1$ we get a
   constant $\gamma$ which does not depend on the chosen divisor
   $D$, while in the remaining cases  the ratio of $a$ and $b$ is
   involved in $\gamma$. This means that an asymptotical behaviour can
   only be examined if the ratio remains unchanged.

%%%%%%%%%%%%%%%%%%%%%%%%%%%%%%%%%%%%%%%%%%%%%%%%%%%%%%%%%%%%%%%%%%%%%%%%%%%%%%%%%%%%%%%%%%%%

   \subsection{Geometrically Ruled Surfaces}\label{subsec:irreducibility:ruled-surfaces}

   Let $\pi:\Sigma=\P(\ke)\rightarrow C$ be a geometrically ruled
   surface with normalised bundle $\ke$ (in the 
   sense of \cite{Har77} V.2.8.1). The N\'eron--Severi group of
   $\Sigma$ is $\NS(\Sigma) = C_0\Z\oplus F\Z$ 
%   \begin{displaymath}
%     \NS(\Sigma) = C_0\Z\oplus F\Z,
%   \end{displaymath}
   with intersection matrix $\left(\begin{smallmatrix}-e & 1 \\ 1 & 0\end{smallmatrix}\right)$
%   \begin{displaymath}
%     \left(\!\begin{array}{rc}-e & 1 \\ 1 & 0\end{array}\right),
%   \end{displaymath}
   where $F\cong\PC^1$ is a fibre of $\pi$, $C_0$ a section of $\pi$
   with $\ko_\Sigma(C_0)\cong\ko_{\P(\ke)}(1)$, $g=g(C)$ the genus of
   $C$,  $\mathfrak{e}=\Lambda^2\ke$ and
   $e=-\deg(\mathfrak{e})\geq -g$. 
   For the canonical divisor we have $K_\Sigma \sim_a -2C_0+ (2g-2-e)\cdot F$.
%   \begin{displaymath}
%     K_\Sigma \sim_a -2C_0+ (2g-2-e)\cdot F,
%   \end{displaymath}
%   where $g=g(C)$ is the genus of the base curve $C$.

   \begin{theorem}\label{thm:irred-ruled-surfaces}
     Let $\pi:\Sigma\rightarrow C$ be a geometrically ruled surface
     with  $e\leq 0$.
     Let $D=aC_0+bF\in\Div(\Sigma)$ with $a\geq 2$, $b>2g-2+\frac{ae}{2}$, and
     if $g=0$ then $b\geq 2$. Let
     $\ks_1,\ldots,\ks_r$ be pairwise distinct
     topological or analytical singularity types and
     $k_1,\ldots,k_r\in\N\setminus\{0\}$. 

     Suppose that 
     \begin{equation}
       \label{eq:irred-ruled-surfaces:1}
       \sum\limits_{i=1}^rk_i\deg\big(X(\ks_i)\big)^2
       \;<\;\gamma\cdot(D-K_\Sigma)^2, 
     \end{equation} 
     where $\gamma$ may be taken from the following table with
     $\alpha=\frac{a+2}{b+2-2g-\tfrac{ae}{2}}>0$. 
     \begin{center}
       \begin{tabular}{|r|r|c|}
         \hline
         \multicolumn{1}{|c}{$g$}&\multicolumn{1}{|c|}{$e$}&$\gamma$\\\hline\hline
          $0$ & $0$ &  $\tfrac{1}{24}$ \\
          $1$ & $0$ & $\tfrac{1}{\max\{24,2\alpha\}}$ \\
          $1$ & $-1$ & $\tfrac{1}{\max\Big\{\min\big\{30+\tfrac{16}{\alpha}+4\alpha,40+9\alpha\big\},
           \tfrac{13}{2}\alpha\Big\}}$\\ 
          $\geq 2 $&$ 0 $&$  \tfrac{1}{\max\{24+16g, 4g\alpha\}}$\\ 
          $\geq 2 $&$ <0 $& 
          $ \tfrac{1}{\max\Big\{\min\big\{24+16g-9e\alpha,18+16g-9e\alpha
           -\tfrac{16}{e\alpha}\big\},4g\alpha-9e\alpha\Big\}}$ \\\hline
       \end{tabular}
%       \tablefirsthead{\hline
%         \multicolumn{1}{|c}{$g$}&\multicolumn{1}{|c|}{$e$}&$\gamma$&$\gamma$, if $\alpha=1$\\\hline\hline}
%       \tablehead{\hline
%         \multicolumn{1}{|c}{$g$}&\multicolumn{1}{|c|}{$e$}&$\gamma$&$\gamma$, if $\alpha=1$\\\hline\hline}
%       \tabletail{\hline}
%       \tablelasttail{\\\hline}
%       \begin{supertabular}{|r|r|c|c|}         
%          $0$ & $0$ &  $\tfrac{1}{24}$ & $\frac{1}{24}$ \\
%          $1$ & $0$ & $\tfrac{1}{\max\{24,2\alpha\}}$ & $\frac{1}{24}$\\
%          $1$ & $-1$ & $\tfrac{1}{\max\Big\{\min\big\{30+\tfrac{16}{\alpha}+4\alpha,40+9\alpha\big\},
%           \tfrac{13}{2}\alpha\Big\}}$&$\frac{1}{49} $\\ 
%          $\geq 2 $&$ 0 $&$  \tfrac{1}{\max\{24+16g, 4g\alpha\}}$&$\frac{1}{24+16g} $\\ 
%          $\geq 2 $&$ <0 $& 
%          $ \tfrac{1}{\max\Big\{\min\big\{24+16g-9e\alpha,18+16g-9e\alpha
%           -\tfrac{16}{e\alpha}\big\},4g\alpha-9e\alpha\Big\}}$ &
%       \end{supertabular}
     \end{center}

     Then either
     $V_{|D|}^{irr}(k_1\ks_1,\ldots,k_r\ks_r)$ is empty
     or it is irreducible of the expected dimension.
     \hfill $\Box$
   \end{theorem}

   Once more, only in the case $g=0$, i.\ e.\ when $\Sigma\cong\PC^1\times\PC^1$, we
   are in the lucky situation that the 
   constant $\gamma$ does not at all depend on the chosen divisor
   $D$, whereas in the case $g\geq 1$ the ratio of $a$ and $b$ is
   involved in $\gamma$. This means that an asymptotical behaviour can
   only be examined if the ratio remains unchanged.

   If $\Sigma$ is a product $C\times\PC^1$ the constant $\gamma$ here is
   the same as in Section 
   \ref{subsec:irreducibility:product-surfaces}. 

   In \cite{Ran89} and in \cite{GLS98b} the case of nodal curves on the Hirzebruch
   surface $\F_1$ is treated, since this is just $\PC^2$ blown up at
   one point. $\F_1$ is an example of a geometrically ruled surface
   with invariant $e=1>0$, a case which we so far cannot treat with our
   methods, due to the section with self-intersection $-1$. However, 
   it seems to be
   possible to extend the methods of \cite{GLS98b} to the situation of
   arbitrary ruled surfaces with positive invariant $e$ -- at least if
   we restrict to singularities which are not too bad.

%%%%%%%%%%%%%%%%%%%%%%%%%%%%%%%%%%%%%%%%%%%%%%%%%%%%%%%%%%%%%%%%%%%%%%%%%%%%%%%%%%%%%%%%%%%%

   \subsection{The Proofs}

   Our approach to the problem proceeds along the lines of an unpublished result of
   Greuel, Lossen and Shustin (cf.~\cite{GLS98d}), which is based on
   ideas of Chiantini and Ciliberto (cf.~\cite{CC99}). The basic ideas are
   in some respect similar to the approach used in \cite{GLS00},
   replacing the ``Castelnuovo-function'' arguments by ``Bogomolov
   instability''. 
   
   We first show that the open subscheme $V^{irr,reg}=V_{|D|}^{irr,reg}(k_1\ks_1,\ldots,k_r\ks_r)$ of 
   $V^{irr}=V_{|D|}^{irr}(k_1\ks_1,\ldots,k_r\ks_r)$,
   and hence its 
   closure $\overline{V^{irr,reg}}$ in $V^{irr}$, is always irreducible (cf.~Theorem
   \ref{thm:v-reg}), 
   and then we look for criteria which
   ensure that the complement of $\overline{V^{irr,reg}}$ in $V^{irr}$
   is empty (cf.~Section \ref{sec:technical-lemmata}). For the latter, we consider
   the restriction of the 
   morphism 
     $\Psi:V\rightarrow B$ (cf.~Subsection \ref{subsec:psi})
   to an irreducible component $V^*$ of $V^{irr}$ not contained in
   $\overline{V^{irr,reg}}$. From the fact that the dimension of $V^*$ is
   at least the expected dimension $\dim\big(V^{irr,reg}\big)$, we
   deduce that the codimension of $B^*=\Psi\big(V^*\big)$ in $B$ is 
   at most $h^1\big(\Sigma,\kj_{X(C)/\Sigma}(D)\big)$, where $C\in
   V^*$ (cf.~Lemma \ref{lem:irred-E}). It thus suffices to find
   conditions which contradict this 
   inequality, that is, we have to get our hands on $\codim_B(B^*)$. 
   However, on the surfaces which we consider the non-vanishing of
   $h^1\big(\Sigma,\kj_{X(C)/\Sigma}(D)\big)$ means in some sense that
   the zero-dimensional scheme $X(C)$ is in special position. We may
   thus hope to realise large parts $X_i^0$ of $X(C)$ on curves
   $\Delta_i$ of ``small degree'' ($i=1,\ldots,m$), which would impose
   at least $\# X_i^0-\dim|\Delta_i|_l$ conditions on $X(C)$, giving
   rise to a lower bound $\sum_{i=1}^m \# X_i^0-\dim|\Delta_i|_l$ for
   $\codim_B(B^*)$. The $X_i^0$'s
   and the $\Delta_i$'s are found in Lemma \ref{lem:irred-A} with
   the aid of certain Bogomolov\index{Bogomolov} unstable rank-two bundles. It thus finally
   remains (cf.~Lemma \ref{lem:irred-C}, \ref{lem:irred-D} and
   \ref{lem:irred-D*}) to give conditions which imply
   \begin{displaymath}
     \sum_{i=1}^m \#
     X_i^0-\dim|\Delta_i|_l>h^1\big(\Sigma,\kj_{X(C)/\Sigma}(D)\big). 
   \end{displaymath}
   These considerations lead to the following proofs.

   \begin{proof}[Proof of Theorem \ref{thm:irred-p3}]
     We may assume that
     $V^{irr}$
     is non-empty.
     By Theorem \ref{thm:v-reg} it suffices to show that
     $V^{irr}=\overline{V^{irr,reg}}$.

     Suppose the contrary, i.e., there is an irreducible curve
     $C_0\in V^{irr}\setminus\overline{V^{irr,reg}}$, in particular
     $h^1\big(\Sigma,\kj_{X_0/\Sigma}(D)\big)>0$ for $X_0=X(C_0)$. Since
     $\deg(X_0)=\sum_{i=1}^rk_i\deg\big(X(\ks_i)\big)$ and 
     $\sum_{z\in\Sigma}\big(\deg(X_{0,z})\big)^2=
     \sum\limits_{i=1}^rk_i\deg(X\big(\ks_i)\big)^2$ the assumptions (0)-(3) of
     Lemma \ref{lem:irred-A} and (4) of Lemma \ref{lem:irred-C} are
     fulfilled. Thus Lemma \ref{lem:irred-C} implies that $C_0$
     satisfies Condition \eqref{eq:irred-E:1} in Lemma
     \ref{lem:irred-E}, which it cannot
     satisfy by the same Lemma. Thus we have derived a
     contradiction.
   \end{proof}

   \begin{proof}[Proof of Theorem \ref{thm:irred-products-of-curves}]
     The assumptions on $a$ and $b$ ensure that $D-K_\Sigma$ is big
     and nef and that $D+K_\Sigma$ is nef. Thus,
     once we know that \eqref{eq:irred-products-of-curves:1} implies Condition (3)
     in Lemma \ref{lem:irred-A} we can do the same proof as in Theorem
     \ref{thm:irred-p3}, just replacing Lemma~ \ref{lem:irred-C} by
     Lemma~ \ref{lem:irred-D}.

     For Condition (3) we note that
     \begin{displaymath}
         \sum\limits_{i=1}^rk_i\deg\big(X(\ks_i)\big)
         \leq\sum\limits_{i=1}^rk_i\cdot\Big(\deg\big(X(\ks_i)\big)\Big)^2 
         \leq\tfrac{1}{24}\cdot(D-K_\Sigma)^2
         <\tfrac{1}{4}\cdot(D-K_\Sigma)^2.
     \end{displaymath}
   \end{proof}

   \begin{proof}[Proof of Theorem \ref{thm:irred-ruled-surfaces}]
     The proof is identical to that of Theorem
     \ref{thm:irred-products-of-curves}, just replacing 
     Lemma \ref{lem:irred-D} by Lemma \ref{lem:irred-D*}.
   \end{proof}

%%%%%%%%%%%%%%%%%%%%%%%%%%%%%%%%%%%%%%%%%%%%%%%%%%%%%%%%%%%%%%%%%%%%%%%%%%%%%%%%%%%%%%%%

   \subsection{Some Remarks}

   What are the obstructions to our approach?

   First, the Bogomolov\index{Bogomolov} instability does not give  much information about
   the curves $\Delta_i$ apart from their existence and the fact that
   they are in some sense ``small'' compared with the divisor $D$. We
   are thus  bound to the study of surfaces where we
   have a good knowledge of the dimension of arbitrary complete linear
   systems. Second, in order to derive the above inequality many
   nasty calculations are necessary which strongly depend on
   the particular structure of the N\'eron--Severi group of the
   surface, that is, we are restricted to surfaces where the
   N\'eron--Severi group is not too large
   and the intersection pairing is not too hard
   (cf.~Lemma \ref{lem:irred-C}, \ref{lem:irred-D} and
   \ref{lem:irred-D*}). Finally, in order to ensure the Bogomolov\index{Bogomolov}
   instability of the vector bundle considered throughout the proof of
   Lemma \ref{lem:irred-A} we heavily use the fact that the surface
   $\Sigma$ does not contain any curve of negative self-intersection,
   which excludes e.~g.~general Hirzebruch surfaces. 

   If the number of irreducible curves of negative self-intersection
   is not too large, one might overcome this last obstacle with the
   technique used in \cite{GLS98b}. That is, we would have to show
   that under certain additional conditions the singular points of the
   considered curves could be independently moved, in particular, they
   could be moved off the exceptional curves - more precisely, the
   subvariety of $V^{irr}$ of curves whose singular locus does not lie
   on any exceptional curve is dense in $V^{irr}$. For this one
   basically just needs criteria for the existence of 
   ``small'' curves realising a zero-dimensional scheme slightly bigger
   than the  equisingularity scheme (respectively the equianalytical
   singularity scheme) of the members in
   $V^{irr}$. E.~g.~in the case of curves with $r$ nodes, that means the
   existence of curves passing through $r$ arbitrary points and having
   multiplicity two in one of them. 
   
   In Section \ref{sec:v-irr-reg} we not only prove that $V^{irr,reg}$
   is irreducible, but also that this indeed remains true if we drop
   the requirement that the curves should be irreducible,
   i.~e.~we show that $V^{reg}$ is irreducible. However, unfortunately our
   approach does not give conditions for the emptiness of the
   complement of $\overline{V^{reg}}$, and thus we cannot say anything
   about the irreducibility of the variety of possibly reducible
   curves in $|D|_l$ with prescribed singularities.
   The reason for this is that in the proof of
   Lemma \ref{lem:irred-A}  we use the
   Theorem of B\'ezout to estimate $D.\Delta_i$. 
%   Since
%   $\Delta_i$ may be about ``half'' of $D$, we need an irreducible
%   curve in $|D|_l$ to be sure that at least for some curve in $|D|_l$ the
%   curve $\Delta_i$ is not a component.
   
%%%%%%%%%%%%%%%%%%%%%%%%%%%%%%%%%%%%%%%%%%%%%%%%%%%%%%%%%%%%%%%%%%%%%%%%%%%%%%%%%%%%%%
   
   \section{$V^{irr,reg}$ is irreducible}\label{sec:v-irr-reg}

   We now show that $V^{irr,reg}$ is always
   irreducible. We do this by showing that under $\Psi:V\rightarrow B$
   every irreducible component of $V^{irr,reg}$ is smooth and maps dominant to the
   irreducible variety $B$ with irreducible fibres.

   \begin{theorem}\label{thm:v-reg}
     Let $D\in\Div(\Sigma)$, $\ks_1,\ldots,\ks_r$ be pairwise distinct
     topological or analytical singularity types and
     $k_1,\ldots,k_r\in\N\setminus\{0\}$.

     If
     $V_{|D|}^{irr,reg}(k_1\ks_1,\ldots,k_r\ks_r)$ is
     non-empty, then it is a T-smooth, irreducible, open subset of
     $V^{irr}_{|D|}(k_1\ks_1,\ldots,k_r\ks_r)$
     of dimension\tom{\footnote{Here $\deg\big(X^*(\ks_i)\big)$ means
       $\deg\big(X^{es}(\ks_i)\big)$ respectively
       $\deg\big(X^{ea}(\ks_i)\big)$.}}
     $\dim |D|_l-\sum_{i=1}^rk_i\deg\big(X^*(\ks_i)\big)$.
   \end{theorem}
   \begin{proof}
     Since $V_{|D|}^{irr,reg}(k_1\ks_1,\ldots,k_r\ks_r)$ is an open
     subset of $V_{|D|}^{reg}(k_1\ks_1,\ldots,k_r\ks_r)=V^{reg}$, it
     suffices to show the claim for $V^{reg}$.

     Let us consider
     the following maps from Subsection \ref{subsec:psi}
     \begin{displaymath}
         \xymatrix{
           \Psi=\Psi_D(k_1\ks_1,\ldots,k_r\ks_r):
           V=V_{|D|}(k_1\ks_1,\ldots,k_r\ks_r)\ar[r]&
           B(k_1\ks_1,\ldots,k_r\ks_r)
           }       
     \end{displaymath}
     and 
     \begin{displaymath}
         \xymatrix{
           \psi=\psi(k_1\ks_1,\ldots,k_r\ks_r):
           B(k_1\ks_1,\ldots,k_r\ks_r)\ar[r]&
           \Hilb_\Sigma^n.
           }       
     \end{displaymath}

     \begin{varthm-roman}[Step 1]
       Every irreducible component $V^*$ of $V^{reg}$ is T-smooth of
       dimension $\dim |D|_l-\sum_{i=1}^rk_i\deg\big(X^*(\ks_i)\big)$. 
     \end{varthm-roman}          
     By \cite{Los98} Proposition 2.1 (c2) $V^*$ is T-smooth at any
     $C\in V^*$
     of dimension $\dim |D|_l-\deg\big(X^*(C)\big)$, since
     $h^1\big(\Sigma,\kj_{X^*/\Sigma}(D)\big)=0$\tom{ according
     to Remark~ \ref{rem:zero-dim-scheme}}. Note that
     $\deg\big(X^*(C)\big)=\sum_{i=1}^rk_i\deg\big(X^*(\ks_i)\big)$ only depends on
     $k_1\ks_1,\ldots,k_r\ks_r$ (cf.~Subsection
     \ref{subsec:schemes}). 

     \begin{varthm-roman}[Step 2]
       $V^{reg}$ is open in $V$. 
     \end{varthm-roman}
     Let $C\in V^{reg}$, then
     $h^1\big(\Sigma,\kj_{X(C)/\Sigma}(D)\big)=0$.
     Thus by semicontinuity  there exists an open, dense
     neighbourhood $U$ of $X(C)$ in
     $\Hilb_\Sigma^n$ such that
     $h^1\big(\Sigma,\kj_{Y/\Sigma}(D)\big)=0$ for all $Y\in U$. But
     then $\Psi^{-1}\big(\psi^{-1}(U)\big)\subseteq V^{reg}$ is an
     open neighbourhood of $C$ in
     $V$, and hence
     $V^{reg}$ is open in $V$.

     \begin{varthm-roman}[Step 3]
       $\Psi$
       restricted to any irreducible component $V^*$ of $V^{reg}$ is dominant.
     \end{varthm-roman}
     Let $V^*$ be an irreducible component of $V^{reg}$ and let $C\in
     V^*$. Since $\Psi^{-1}\big(\Psi(C)\big)$ is an open, dense
     subset of $\big|\kj_{X(C)/\Sigma}(D)\big|_l$
     and since $h^1\big(\Sigma,\kj_{X(C)/\Sigma}(D)\big)=0$, we have
     \begin{displaymath}
       \dim\Psi^{-1}\big(\Psi(C)\big)=h^0\big(\Sigma,\kj_{X(C)/\Sigma}(D)\big)-1
       =\dim|D|_l-\deg\big(X(C)\big).
     \end{displaymath}
     By Step 1 we know the dimension of $V^*$ and by 
     Subsection \ref{subsec:psi} we also know the dimension of
     $B$. Thus we conclude
     \begin{displaymath}
       \begin{array}{rcl}
         \dim\Psi\big(V^*\big)&=&\dim V^*-\dim\Psi^{-1}\big(\Psi(C)\big) \\
         &=&\big(\dim|D|_l-\deg X^*(C)\big)
         -\big(\dim|D|_l-\deg X(C)\big)\\ 
         &=&\deg\big(X(C)\big)-\deg\big(X^*(C)\big)        
         \;\,=\;\,\dim B.
       \end{array}
     \end{displaymath}
     Since $B$ is irreducible $\Psi\big(V^*\big)$ must be
     dense in $B$.

     \begin{varthm-roman}[Step 4]
       $V^{reg}$ is irreducible.
     \end{varthm-roman}
     Let $V^*$ and $V^{**}$ be two irreducible components of
     $V^{reg}$. Then
     $\Psi\big(V^*\big)\cap\Psi\big(V^{**}\big)\not=\emptyset$,
     and thus some fibre $F$ of $\Psi$ intersects both, $V^*$ and
     $V^{**}$. However, the fibre is irreducible and by Step 1 both
     $V^*$ and $V^{**}$ are smooth. Thus $F$ must be completely
     contained in $V^*$ and $V^{**}$, which implies that $V^*=V^{**}$, since both are
     smooth of the same dimension. Thus
     $V^{reg}$ is irreducible.
   \end{proof}

%%%%%%%%%%%%%%%%%%%%%%%%%%%%%%%%%%%%%%%%%%%%%%%%%%%%%%%%%%%%%%%%%%%%%%%%%%%%%%%%%

   \section{The Technical Details}\label{sec:technical-lemmata}

   The following lemma is the heart of the proof. Given a curve
   $C\in|D|_l$, whose (analytical) singularity scheme $X_0=X(C)$ 
   is \emph{special with respect to $D$} in the
   sense that $h^1\big(\Sigma,\kj_{X_0/\Sigma}(D)\big)>0$, 
   provides a ``small'' curve $\Delta_1$ through a subscheme $X_1^0$
   of $X_0$, so that we can reduce the problem by replacing $X_0$ and
   $D$ by $X_0:\Delta_1$ and $D-\Delta_1$ respectively. We can of
   course proceed inductively as long as the new zero-dimensional
   scheme is again special with respect to the new divisor.

   In order to find $\Delta_1$ we choose a subscheme
   $X_1^0\subseteq X_0$ which is minimal among those subschemes special with
   respect to $D$. By Grothendieck-Serre duality 
   \begin{displaymath}
     H^1\big(\Sigma,\kj_{X_1^0/\Sigma}(D)\big)\cong
     \Ext^1\big(\kj_{X_1^0/\Sigma}(D-K_\Sigma),\ko_\Sigma\big) 
   \end{displaymath}
   and a non-trivial element of the latter group gives rise to an
   extension
   \begin{displaymath}
     0\rightarrow \ko_\Sigma\rightarrow E_1\rightarrow
     \kj_{X_1^0/\Sigma}(D-K_\Sigma) \rightarrow 0.
   \end{displaymath}
   We then show that the rank-two bundle  $E_1$ is Bogomolov\index{Bogomolov} unstable and deduce the
   existence of a divisor $\Delta_1^0$ such that
   \begin{displaymath}
     H^0\Big(\Sigma,\kj_{X_1^0/\Sigma}\big(D-K_\Sigma-\Delta_1^0\big)\Big)\not=0,
   \end{displaymath}
   that is, we find a curve $\Delta_1\in
   \big|\kj_{X_1^0/\Sigma}\big(D-K_\Sigma-\Delta_1^0\big)\big|_l$. 

   \begin{lemma}\label{lem:irred-A}
     Let $\Sigma$ be a surface such that  any curve $C\subset\Sigma$
     is nef (*).
   
     Let $D\in\Div(\Sigma)$ and $X_0\subset\Sigma$  a
     zero-dimensional scheme satisfying
     \begin{myenumerate}
     \item[\rm(0)] $D-K_\Sigma$ is big and nef,  and $D+K_\Sigma$ is nef,
     \item[\rm(1)] $\exists\;C_0\in|D|_l\;\mbox{ irreducible}:\;X_0\subset C_0$,
     \item[\rm(2)] $h^1\big(\Sigma,\kj_{X_0/\Sigma}(D)\big)>0$, and 
     \item[\rm(3)] $\deg(X_0)<\beta\cdot (D-K_\Sigma)^2$ for some $0<\beta\leq
       \frac{1}{4}$.
     \end{myenumerate}
     Then there exist curves $\Delta_1,\ldots,\Delta_m\subset\Sigma$ and
     zero-dimensional locally complete intersections $X_i^0\subseteq
     X_{i-1}\cap\Delta_i$ for $i=1,\ldots,m$, where $X_i=X_{i-1}:\Delta_i$ for
     $i=1,\ldots,m$,  such that 
     \begin{enumerate}
     \item $h^1\Big(\Sigma,\kj_{X_m/\Sigma}\big(D-\sum_{i=1}^m\Delta_i\big)\Big)=0$,
     \end{enumerate}
     and for $i=1,\ldots,m$ 
     \begin{enumerate}\stepcounter{enumi}
     \item $h^1\Big(\Sigma,\kj_{X_i^0/\Sigma}\big(D-\sum_{k=1}^{i-1}\Delta_k\big)\Big)=1$ 
     \item
       $D.\Delta_i\geq\deg(X_{i-1}\cap\Delta_i)\geq\deg\big(X_i^0\big)
       \geq\big(D-K_\Sigma-\sum_{k=1}^i\Delta_k\big).\Delta_i\geq\Delta_i^2\geq 0$ 
     \item $\big(D-K_\Sigma-\sum_{k=1}^i\Delta_k-\Delta_i\big)^2>0$,
     \item  $\big(D-K_\Sigma-\sum_{k=1}^i\Delta_k-\Delta_i\big).H>0$\;\;
       for all $H\in\Div(\Sigma)$ ample, and
     \item %$D-\sum_{k=1}^i\Delta_k$ is effective and
       $D-K_\Sigma-\sum_{k=1}^i\Delta_k$ is big and nef.
     \end{enumerate}
     Moreover, it follows
     \begin{equation}
       \label{eq:irred-A:1}
       0\leq \tfrac{1}{4}(D-K_\Sigma)^2-\sum_{i=1}^m\deg\big(X_i^0\big)
       \leq \left(\tfrac{1}{2}(D-K_\Sigma)-\sum_{i=1}^m\Delta_i\right)^2.
     \end{equation}
   \end{lemma}

   \begin{proof}
     We are going to find the schemes $\Delta_i$ and $X_i^0$ recursively. Let
     us therefore suppose that we have already found
     $\Delta_1,\ldots,\Delta_{i-1}$ and $X_1^0,\ldots,X_{i-1}^0$
     satisfying (b)-(f), and suppose that still
     $h^1\Big(\Sigma,\kj_{X_{i-1}/\Sigma}\big(D-\sum_{i=1}^{i-1}\Delta_i\big)\Big)>0$. 
     
     \vspace*{-0.1cm}
     We choose $X_i^0\subseteq X_{i-1}$ minimal such
     that $h^1\Big(\Sigma,\kj_{X_i^0/\Sigma}\big(D-\sum_{k=1}^{i-1}\Delta_k\big)\Big)>0$. 
     
     \vspace*{-0.2cm}
     \begin{varthm-roman}[Step 1]
       $h^1\Big(\Sigma,\kj_{X_i^0/\Sigma}\big(D-\sum_{k=1}^{i-1}\Delta_k\big)\Big)=1$, 
       i.~e.~(b) is fulfilled.
     \end{varthm-roman}
     \vspace*{-0.1cm}
     Suppose it was strictly larger than one. By (0) respectively (f),
     and by  the
     Ka\-wa\-ma\-ta--Viehweg Vanishing Theorem we have
     $h^1\Big(\Sigma,\ko_\Sigma\big(D-\sum_{k=1}^{i-1}\Delta_k\big)\Big)=0$. 

     Thus $X_i^0$ cannot be empty, that is $\deg\big(X_i^0\big)\geq
     1$ and we may choose a subscheme $Y\subset X_i^0$ of degree
     $\deg(Y)=\deg\big(X_i^0\big)-1$. The inclusion
     $\kj_{X_i^0}\hookrightarrow \kj_Y$ implies
     $h^0\Big(\Sigma,\kj_{X_i^0/\Sigma}\big(D-\sum_{k=1}^{i-1}\Delta_k\big)\Big)
     \leq
     h^0\Big(\Sigma,\kj_{Y/\Sigma}\big(D-\sum_{k=1}^{i-1}\Delta_k\big)\Big)$
     and the structure sequences of $Y$ and $X_i^0$ thus lead to
     \begin{displaymath}
       h^1\Big(\Sigma,\kj_{Y/\Sigma}\big(D-\mbox{$\sum_{k=1}^{i-1}$}\Delta_k\big)\Big)
       \geq
       h^1\Big(\Sigma,\kj_{X_i^0/\Sigma}\big(D-\mbox{$\sum_{k=1}^{i-1}$}\Delta_k\big)\Big)-1
       >0
%       \as{1.4}
%       \begin{array}{ll}
%         \multicolumn{2}{l}{
%           h^1\Big(\Sigma,\kj_{Y/\Sigma}\big(D-\sum_{k=1}^{i-1}\Delta_k\big)\Big)
%           }\\
%         =&h^0\Big(\Sigma,\kj_{Y/\Sigma}\big(D-\sum_{k=1}^{i-1}\Delta_k\big)\Big)
%         -h^0\Big(\Sigma,\ko_\Sigma\big(D-\sum_{k=1}^{i-1}\Delta_k\big)\Big)
%         +\deg(Y)\\
%         \geq&
%         h^0\Big(\Sigma,\kj_{X_i^0/\Sigma}\big(D-\sum_{k=1}^{i-1}\Delta_k\big)\Big)
%         -h^0\Big(\Sigma,\ko_\Sigma\big(D-\sum_{k=1}^{i-1}\Delta_k\big)\Big)
%         +\deg\big(X_i^0\big)-1\\
%         =&h^1\Big(\Sigma,\kj_{X_i^0/\Sigma}\big(D-\sum_{k=1}^{i-1}\Delta_k\big)\Big)-1
%         >0
%       \end{array}
%       \as{\asf}
     \end{displaymath}
     contradicting the minimality of $X_i^0$.

     \begin{varthm-roman}[Step 2]
       $\deg\big(X_i^0\big)\leq \deg(X_0)-\sum_{k=1}^{i-1}\deg(X_{k-1}\cap\Delta_k)$.
     \end{varthm-roman}
     The case $i=1$ follows from the
     fact that $X_1^0\subseteq X_0$, and for $i>1$
     the inclusion $X_i^0\subseteq X_{i-1}=X_{i-2}:\Delta_{i-1}$
     implies
     \begin{displaymath}
       \deg\big(X_i^0\big)\leq \deg(X_{i-2}:\Delta_{i-1})=\deg(X_{i-2})-\deg(X_{i-2}\cap\Delta_{i-1}).
     \end{displaymath}
     It thus suffices to show, that 
     \begin{displaymath}
       \deg(X_{i-2})-\deg(X_{i-2}\cap\Delta_{i-1})
       =
       \deg(X_0)-\mbox{$\sum\nolimits_{k=1}^{i-1}$}\deg(X_{k-1}\cap\Delta_k).
     \end{displaymath}
     If $i=2$, there is nothing to show. Otherwise
     $X_{i-2}=X_{i-3}:\Delta_{i-2}$ implies
     \begin{displaymath}
       \begin{array}{ll}
         \multicolumn{2}{l}{\deg(X_{i-2})-\deg(X_{i-2}\cap\Delta_{i-1})}\\
         =& 
         \deg(X_{i-3}:\Delta_{i-2})-\deg(X_{i-2}\cap\Delta_{i-1})\\
         =&
         \deg(X_{i-3})-\deg(X_{i-3}\cap\Delta_{i-2})-\deg(X_{i-2}\cap\Delta_{i-1})
       \end{array}
     \end{displaymath}
     and we are done by induction.

     \begin{varthm-roman}[Step 3]
       There exists a ``suitable'' %, i.~e.~satisfying \eqref{eq:irred-A:2} and \eqref{eq:irred-A:3}, 
       locally free rank-two vector bundle $E_i$.
     \end{varthm-roman}
     By the Grothendieck-Serre duality \tom{(cf.~\cite{Har77}
       III.7.6) }
     we have
     \begin{displaymath}
       0\not=H^1\Big(\Sigma,\kj_{X_i^0/\Sigma}\big(D-\mbox{$\sum_{k=1}^{i-1}$}\Delta_k\big)\Big)
       \cong\Ext^1\Big(\kj_{X_i^0/\Sigma}\big(D-K_\Sigma-\mbox{$\sum_{k=1}^{i-1}\Delta_k$}\big),\ko_\Sigma\Big).
     \end{displaymath}
     That is, there exists an extension  \tom{(cf.~\cite{Har77}
       Ex.~III.6.1) }
     \begin{equation}\label{eq:irred-A:2}
       0\rightarrow \ko_\Sigma\rightarrow
       E_i\rightarrow\kj_{X_i^0/\Sigma}\left(D-K_\Sigma-\mbox{$\sum\nolimits_{k=1}^{i-1}\Delta_k$}\right) 
       \rightarrow 0.
     \end{equation}
     The minimality of $X_i^0$ implies that $E_i$ is locally free
     \tom{(cf.~\cite{Laz97} Proposition~3.9)} and hence that $X_i^0$ is a
     locally complete intersection (cf.~\cite{Laz97}\tom{ p.~175}). 
     Moreover, we have \tom{(cf.~\cite{Laz97} Exercise 4.3)}
     \begin{equation}\label{eq:irred-A:3}
       c_1(E_i)=D-K_\Sigma-\sum_{k=1}^{i-1}\Delta_k\mbox{\;\;\;and\;\;\;}
       c_2(E_i)=\deg\big(X_i^0\big).
     \end{equation}
     
     \begin{varthm-roman}[Step 4]
       $E_i$ is Bogomolov\index{Bogomolov} unstable.
     \end{varthm-roman}
     According to the Theorem of Bogomolov\index{Bogomolov} we only have to show
     $c_1(E_i)^2>4c_2(E_i)$ (cf.~\cite{Bog79} or \cite{Laz97} Theorem 4.2).
%     In order to show this we note that
%     \begin{equation}
%       \label{eq:irred-A:4}
%       \as{1.4}
%       \begin{array}{ll}
%         \multicolumn{2}{l}{-2\sum_{k=1}^{i-1}\Delta_k.\big(D-K_\Sigma-\sum_{j=1}^k\Delta_j\big)}\\
%         =&
%         -2\left(\sum_{k=1}^{i-1}\Delta_k\right).(D-K_\Sigma)
%         +
%         2\sum_{k=1}^{i-1}\sum_{j=1}^k \Delta_k.\Delta_j\\
%         =&
%         -2\left(\sum_{k=1}^{i-1}\Delta_k\right).(D-K_\Sigma)
%         +
%         \sum_{k=1}^{i-1}\Delta_k^2
%         +
%         \left(\sum_{k=1}^{i-1}\Delta_k\right)^2.
%       \end{array}
%     \end{equation}
     Since
     $(4\beta-1)\cdot(D-K_\Sigma)^2\leq 0$ by (3) and since $\Delta_k^2\geq
     0$ by (*) we deduce:
     \begin{displaymath}
       \as{1.6}
       \begin{array}{rl}
         4 c_2(E_i)&=\;\, 4\deg\big(X_i^0\big)\;\,\leq_\expl{Step 2}\;\,
         4\deg(X_0)-4\sum_{k=1}^{i-1}\deg(X_{k-1}\cap\Delta_k)\\
         &<_\expl{(3)/(c)}
         4\beta
         (D-K_\Sigma)^2-2\sum_{k=1}^{i-1}\Delta_k.\big(D-K_\Sigma-\sum_{j=1}^k\Delta_j\big) 
         - 2\sum_{k=1}^{i-1}\Delta_k^2\\
%         &=_\expl{\eqref{eq:irred-A:4}}
%         4\beta
%         (D-K_\Sigma)^2-2\left(\sum_{k=1}^{i-1}\Delta_k\right).(D-K_\Sigma)
%         +
%         \left(\sum_{k=1}^{i-1}\Delta_k\right)^2
%         -\sum_{k=1}^{i-1}\Delta_k^2\\
         &=\;\,
         \left(D-K_\Sigma-\sum_{k=1}^{i-1}\Delta_k\right)^2+(4\beta-1)\cdot(D-K_\Sigma)^2
         -\sum_{k=1}^{i-1}\Delta_k^2\\
         &\leq\;\,\left(D-K_\Sigma-\sum_{k=1}^{i-1}\Delta_k\right)^2
         =\;\,
         c_1(E_i)^2.         
       \end{array}
       \as{\asf}
     \end{displaymath}

     \begin{varthm-roman}[Step 5]
       Find $\Delta_i$.
     \end{varthm-roman}
     Since $E_i$ is Bogomolov\index{Bogomolov} unstable there is a $0$-dim.\
     scheme $Z_i\subset\Sigma$ and a 
     $\Delta_i^0\in\Div(\Sigma)$: %s.\ t.\
     \begin{equation}
       \label{eq:irred-A:5}
       0\rightarrow\ko_\Sigma\big(\Delta_i^0\big)\rightarrow
       E_i\rightarrow\kj_{Z_i/\Sigma}\left(D-K_\Sigma-\mbox{$\sum\nolimits_{k=1}^{i-1}\Delta_k$}-\Delta_i^0\right)
       \rightarrow 0
     \end{equation}
     is exact \tom{(cf.~\cite{Laz97} Theorem 4.2)} and such that
     \begin{myenumerate}
     \item[\rm(d')]
       $\big(2\Delta_i^0-D+K_\Sigma+\sum_{k=1}^{i-1}\Delta_k\big)^2\geq c_1(E_i)^2-4\cdot c_2(E_i)>0$, and
     \item[\rm(e')]
       $\big(2\Delta_i^0-D+K_\Sigma+\sum_{k=1}^{i-1}\Delta_k\big).H>0$\;\;
       for all $H\in\Div(\Sigma)$ ample.
     \end{myenumerate}
%     We note that (e') implies 
%     $h^0\Big(\Sigma,\ko_\Sigma\big(D-K_\Sigma
%     -\mbox{$\sum\nolimits_{k=1}^{i-1}\Delta_k$}-2\Delta_i^0\big)\Big)=0$
%     and the inclusion $\kj_{Z_i/\Sigma}\hookrightarrow\ko_\Sigma$
%     thus gives 
%     \begin{equation}\label{eq:irred-A:6*}
%       h^0\Big(\Sigma,\kj_{Z_i/\Sigma}\big(D-K_\Sigma
%       -\mbox{$\sum\nolimits_{k=1}^{i-1}\Delta_k$}-2\Delta_i^0\big)\Big)=0.
%     \end{equation}
     Tensoring \eqref{eq:irred-A:5} with
     $\ko_\Sigma\big(-\Delta_i^0\big)$ leads to the following exact sequence
     \begin{equation}
       \label{eq:irred-A:6}
       0\rightarrow\ko_\Sigma\rightarrow
       E_i\big(-\Delta_i^0\big)
       \rightarrow\kj_{Z_i/\Sigma}\left(D-K_\Sigma-\mbox{$\sum\nolimits_{k=1}^{i-1}\Delta_k$}-2\Delta_i^0\right)
       \rightarrow 0,
     \end{equation}
     and we deduce %with \eqref{eq:irred-A:6*}
     that $h^0\Big(\Sigma,E_i\big(-\Delta_i^0\big)\Big)\not=0$.%=h^0(\Sigma,\ko_\Sigma)=1$. 

     Now tensoring \eqref{eq:irred-A:2} with $\ko_\Sigma\big(-\Delta_i^0\big)$ leads to
     \begin{equation}\label{eq:irred-A:7}
       0\rightarrow \ko_\Sigma\big(-\Delta_i^0\big)\rightarrow
       E_i\big(-\Delta_i^0\big)
       \rightarrow\kj_{X_i^0/\Sigma}\left(D-K_\Sigma-\mbox{$\sum\nolimits_{k=1}^{i-1}\Delta_k$}-\Delta_i^0\right) 
       \rightarrow 0.
     \end{equation}
     By (e'), and (0) respectively (f)
     \begin{displaymath}
       -\Delta_i^0.H<-\tfrac{1}{2}\big(D-K_\Sigma-\mbox{$\sum_{k=1}^{i-1}\Delta_k$}\big).H\leq 0  
     \end{displaymath}
     for an ample divisor $H$, hence
     $-\Delta_i^0$ cannot be effective, that is
     $H^0\big(\Sigma,-\Delta_i^0\big)=0$. But the long exact
     cohomology sequence of \eqref{eq:irred-A:7} then implies
     \begin{displaymath}
       0\not=H^0\Big(\Sigma,E_i\big(-\Delta_i^0\big)\Big)
       \hookrightarrow 
       H^0\left(\Sigma,\kj_{X_i^0/\Sigma}
         \left(D-K_\Sigma-\mbox{$\sum\nolimits_{k=1}^{i-1}\Delta_k$}-\Delta_i^0\right)\right).
     \end{displaymath}
     In particular 
%     the linear system 
%     $\big|\kj_{X_i^0/\Sigma}\big(D-K_\Sigma-\sum\nolimits_{k=1}^{i-1}\Delta_k-\Delta_i^0\big)\big|_l$
%     is non-empty, and 
     we may choose
     $\Delta_i\in
     \Big|\kj_{X_i^0/\Sigma}\big(D-K_\Sigma-
     \mbox{$\sum\nolimits_{k=1}^{i-1}\Delta_k$}-\Delta_i^0\big)\Big|_l.$ 
%     \begin{displaymath}
%       \Delta_i\in
%       \Big|\kj_{X_i^0/\Sigma}\big(D-K_\Sigma-\mbox{$\sum\nolimits_{k=1}^{i-1}\Delta_k$}-\Delta_i^0\big)\Big|_l.
%     \end{displaymath}
     
     \begin{varthm-roman}[Step 6]
       $\Delta_i$ satisfies (d)-(f).
     \end{varthm-roman}
     We note that by the choice of $\Delta_i$ we have the following equivalences
     \begin{equation}
       \label{eq:irred-A:8}
       \Delta_i^0\sim_l D-K_\Sigma-\mbox{$\sum\nolimits_{k=1}^i\Delta_k$}
     \end{equation}
%     and 
     \begin{equation}
       \label{eq:irred-A:9}
       \Delta_i^0-\Delta_i\sim_l2\Delta_i^0-D+K_\Sigma+\mbox{$\sum\nolimits_{k=1}^{i-1}\Delta_k$}
       \sim_lD-K_\Sigma-\mbox{$\sum\nolimits_{k=1}^i\Delta_k$}-\Delta_i.
     \end{equation}
     Thus (d) and (e) is  a reformulation of (d') and (e').

     Moreover, since $\big(\Delta_i^0-\Delta_i\big).H>0$ for any ample $H$,
     then $\big(\Delta_i^0-\Delta_i\big).H\geq 0$ for any $H$ in
     the closure of the ample cone, in particular
%     \begin{equation}\label{eq:irred-A:9*}
%       \big(\Delta_i^0-\Delta_i\big).H\geq 0\;\;\;\mbox{ for all } H
%       \mbox{ nef}.
%     \end{equation}
%     But then 
     \begin{equation}
       \label{eq:irred-A:10}
       \Delta_i^0.H\geq\Delta_i.H\geq 0\;\;\;\mbox{ for all } H
       \mbox{ nef},
     \end{equation}
     since $\Delta_i$ is effective. And finally, since by assumption
     (*) any effective divisor is nef, we deduce that
     $\Delta_i^0.C\geq 0$ for any curve $C$, that is, $\Delta_i^0$ is
     nef. In view of \eqref{eq:irred-A:8} for (f) it remains to show
     that $\big(\Delta_i^0\big)^2>0$. Taking once more into account
     that $\Delta_i$ is nef by (*) we have by (d'), \eqref{eq:irred-A:9},
%     \eqref{eq:irred-A:9*} 
     and \eqref{eq:irred-A:10}
     \begin{displaymath}
       \big(\Delta_i^0\big)^2=
       \big(\Delta_i^0-\Delta_i\big)^2+\big(\Delta_i^0-\Delta_i\big).\Delta_i
       +\Delta_i^0.\Delta_i>0.
     \end{displaymath}

     \begin{varthm-roman}[Step 7]
       $\Delta_i$ satisfies (c).
     \end{varthm-roman}
     We would like to apply the Theorem of B\'ezout to $C_0$ and
     $\Delta_i$. Thus suppose that the irreducible curve $C_0$ is a
     component of $\Delta_i$ and let $H$ be any ample divisor.

     Applying (d) and the fact that $D+K_\Sigma$ is nef by (0), we derive the
     contradiction
     \begin{displaymath}
       0\leq(\Delta_i-C_0).H<-\frac{1}{2}\cdot\left(D+K_\Sigma+\sum_{k=1}^{i-1}\Delta_k\right).H
       \leq -\frac{1}{2}\cdot(D+K_\Sigma).H\leq 0.
     \end{displaymath}
     Since $X_{i-1}\subseteq X_0\subset C_0$ the Theorem of B\'ezout therefore implies
     \begin{displaymath}
       D.\Delta_i=C_0.\Delta_i\geq\deg(X_{i-1}\cap\Delta_i).
     \end{displaymath}
     By definition $X_i^0\subseteq X_{i-1}$ and
     $X_i^0\subset\Delta_i$, thus 
     \begin{displaymath}
       \deg(X_{i-1}\cap\Delta_i)\geq\deg\big(X_i^0\big).
     \end{displaymath}
     By assumption (*) the curve $\Delta_i$ is nef and thus
     \eqref{eq:irred-A:10} gives 
     \begin{displaymath}
       \big(D-K_\Sigma-\mbox{$\sum_{k=1}^i\Delta_k$}\big).\Delta_i=\Delta_i^0.\Delta_i\geq\Delta_i^2\geq 0.
     \end{displaymath}
     Finally from (d') and by \eqref{eq:irred-A:3} it follows that
     \begin{displaymath}
       \big(\Delta_i^0-\Delta_i\big)^2\geq c_1(E_i)^2-4\cdot c_2(E_i)
       =\big(\Delta_i^0+\Delta_i\big)^2-4\cdot\deg\big(X_i^0\big),
     \end{displaymath}
     and thus $\deg\big(X_i^0\big)\geq \Delta_i^0.\Delta_i$.

     \begin{varthm-roman}[Step 8]
       After a finite number $m$ of steps 
       $h^1\Big(\Sigma,\kj_{X_m/\Sigma}\big(D-\sum_{i=1}^m\Delta_i\big)\Big)=0$.
%       i.~e.~(a) is fulfilled.
     \end{varthm-roman}
     As we have mentioned in Step 1 $\deg\big(X_i^0\big)>0$. This
     ensures that 
     \begin{displaymath}
       \deg(X_i)%=\deg(X_{i-1}:\Delta_i)
       =\deg(X_{i-1})-\deg(X_{i-1}\cap\Delta_i)
       \leq\deg(X_{i-1})-\deg\big(X_i^0\big)<\deg(X_{i-1}),
     \end{displaymath}
     i.~e.~the degree of $X_i$ strictly decreases each time.
     Thus the procedure must stop after a finite number $m$ of steps\tom{,
     which is equivalent to the fact that 
     $h^1\Big(\Sigma,\kj_{X_m/\Sigma}\big(D-\sum_{i=1}^m\Delta_i\big)\Big)=0$}. 

     \begin{varthm-roman}[Step 9]
       It remains to show \eqref{eq:irred-A:1}.
%       $0\leq \tfrac{1}{4}(D-K_\Sigma)^2-\sum_{i=1}^m\deg\big(X_i^0\big)
%       \leq \left(\tfrac{1}{2}(D-K_\Sigma)-\sum_{i=1}^m\Delta_i\right)^2$.
     \end{varthm-roman}
     By assumption (*) the curves $\Delta_i$ are nef, in particular
     $\Delta_i.\Delta_j\geq 0$ for all $i,j$. Thus (c) implies
     \begin{displaymath}
       \as{1.4}
       \begin{array}{rcl}
         \sum_{i=1}^m\deg\big(X_i^0\big)&\geq&
         \sum_{i=1}^m\big(D-K_\Sigma-\sum_{k=1}^i\Delta_k\big).\Delta_i\\
%         &=& (D-K_\Sigma).\sum_{i=1}^m\Delta_i - \sum_{1\leq k\leq
%           i\leq m}\Delta_k.\Delta_i\\
         &=& (D-K_\Sigma).\sum_{i=1}^m\Delta_i - \frac{1}{2}
         \left(\big(\sum_{i=1}^m\Delta_i\big)^2+\sum_{i=1}^m\Delta_i^2\right)\\
         &\geq&
         (D-K_\Sigma).\sum_{i=1}^m\Delta_i - \big(\sum_{i=1}^m\Delta_i\big)^2.
       \end{array}
     \end{displaymath}
     But then, taking condition (3) into account,
     \begin{displaymath}
       \begin{array}{rcl}
         0&\leq&\tfrac{1}{4}(D-K_\Sigma)^2-\deg(X_0)
         \leq
         \tfrac{1}{4}(D-K_\Sigma)^2-\sum_{i=1}^m\deg\big(X_i^0\big)\\
         &\leq&
         \tfrac{1}{4}(D-K_\Sigma)^2-
         (D-K_\Sigma).\sum_{i=1}^m\Delta_i +
         \big(\sum_{i=1}^m\Delta_i\big)^2\\
         &=&
         \big(\tfrac{1}{2}(D-K_\Sigma)-\sum_{i=1}^m\Delta_i\big)^2.
       \end{array}
       \as{\asf}
     \end{displaymath}
   \end{proof}

   It is our overall aim to compare the dimension of a cohomology
   group of the form $H^1\big(\Sigma,\kj_{X_0/\Sigma}(D)\big)$ with
   some invariants of the $X_i^0$ and $\Delta_i$. The following lemma
   will be vital for the necessary estimates.

   \begin{lemma}\label{lem:irred-B}
     Let $D\in\Div(\Sigma)$ and let $X_0\subset\Sigma$ be a
     zero-dimensional scheme such that there exist curves 
     $\Delta_1,\ldots,\Delta_m\subset\Sigma$ and
     zero-dimensional schemes $X_i^0\subseteq
     X_{i-1}$ for $i=1,\ldots,m$, where $X_i=X_{i-1}:\Delta_i$ for
     $i=1,\ldots,m$,  such that (a)-(f) in Lemma \ref{lem:irred-A} are
     fulfilled.

     Then:
     \begin{displaymath}
       \as{1.8}
       \begin{array}{rcl}
         h^1\big(\Sigma,\kj_{X_0/\Sigma}(D)\big)&
         \leq&
         \sum\limits_{i=1}^m
         h^1\Big(\Delta_i,\kj_{X_{i-1}\cap\Delta_i/\Delta_i}
         \big(D-\sum_{k=1}^{i-1}\Delta_k\big)\Big)\\
         &\leq&
         \sum\limits_{i=1}^m
         \Big(1+\deg(X_{i-1}\cap\Delta_i)-\deg\big(X_i^0\big)\Big)\\
         &\leq&
         \sum\limits_{i=1}^m
         \Big(\Delta_i\cdot\big(K_\Sigma+\sum_{k=1}^i\Delta_k\big)+1\Big).
       \end{array}
       \as{\asf}
     \end{displaymath}
   \end{lemma}
   \begin{proof}
     Throughout the proof we use the following notation
     \begin{displaymath}
       \kg_i=\kj_{X_{i-1}\cap\Delta_i/\Delta_i}\left(D-\mbox{$\sum\nolimits_{k=1}^{i-1}$}\Delta_k\right)
       \;\;\;\mbox{ and }\;\;\;
       \kg_i^0=\kj_{X_i^0/\Delta_i}\big(D-\mbox{$\sum_{k=1}^{i-1}$}\Delta_k\big)
     \end{displaymath}
     for  $i=1,\ldots,m$, and      for $i=0,\ldots,m$
     \begin{displaymath}
       \kf_i=\kj_{X_i/\Sigma}\left(D-\mbox{$\sum\nolimits_{k=1}^i$}\Delta_k\right).
     \end{displaymath}
     
     Since $X_{i+1}=X_i:\Delta_{i+1}$ we have the following short exact
     sequence
     \begin{equation}
       \xymatrix@C=1cm{
         0\ar[r]& \kf_{i+1}\ar[r]^{\cdot\Delta_{i+1}}&
         \kf_i\ar[r]&\kg_{i+1}\ar[r]&0
         }
     \end{equation}
      for $i=0,\ldots,m-1$ and the corresponding long exact cohomology sequence
     \begin{equation}\label{eq:irred-B:1}
       \begin{aligned}
       \xymatrix@R=0.4cm@C=0.3cm{
        \hspace*{1.3cm}0\ar[r]& H^0(\Sigma,\kf_{i+1})\ar[r]&
         H^0(\Sigma,\kf_i)\ar[r]& H^0(\Sigma,\kg_{i+1})\ar[r]&
         H^1(\Sigma,\kf_{i+1})\ar[d]\\
         0=H^2(\Sigma,\kg_{i+1})
         &H^2(\Sigma,\kf_i)\ar[l]
         &H^2(\Sigma,\kf_{i+1})\ar[l]
         & H^1(\Sigma,\kg_{i+1})\ar[l]
         &H^1(\Sigma,\kf_i)\ar[l]
         }
       \end{aligned}
     \end{equation}

     \begin{varthm-roman}[Step 1]
       $h^1(\Sigma,\kf_i)\leq\sum_{j=i+1}^m h^1(\Sigma,\kg_j)$\;\;
       for $i=0,\ldots,m-1$. 
%       In particular, $h^1\big(\Sigma,\kj_{X_0/\Sigma}(D)\big)
%         \leq
%         \sum_{i=1}^m
%         h^1\Big(\Delta_i,\kj_{X_{i-1}\cap\Delta_i/\Delta_i})$.
     \end{varthm-roman}
     We prove the claim by descending induction on $i$. 
     From \eqref{eq:irred-B:1} we deduce 
     \begin{displaymath}
       \xymatrix@C=0.45cm{
         0=H^1(\Sigma,\kf_m)\ar[r]&
         H^1(\Sigma,\kf_{m-1})\ar[r]& H^1(\Sigma,\kg_m),
         }
     \end{displaymath}
     which implies $h^1(\Sigma,\kf_{m-1})\leq h^1(\Sigma,\kg_m)$ and
     thus proves the case $i=m-1$.
     
     We may
     therefore assume that $1\leq i\leq m-2$. Once more from
     \eqref{eq:irred-B:1}  we deduce
     \begin{displaymath}
       a=h^0(\Sigma,\kf_{i+1})-h^0(\Sigma,\kf_i)+h^0(\Sigma,\kg_{i+1})\geq0,
     \end{displaymath}
     and
     \begin{displaymath}
       b=h^2(\Sigma,\kf_{i+1})-h^2(\Sigma,\kf_i)\geq0, 
     \end{displaymath}
     and finally
     \begin{displaymath}
       \as{1.4}
       \begin{array}{rl}
         h^1(\Sigma,\kf_i)&=\;\,h^1(\Sigma,\kg_{i+1})+h^1(\Sigma,\kf_{i+1})-a-b\;\,
         \leq\;\,h^1(\Sigma,\kg_{i+1})+h^1(\Sigma,\kf_{i+1})\\
         &\leq_\expl{Ind.}\;\,h^1(\Sigma,\kg_{i+1})+\sum_{j=i+2}^m
         h^1(\Sigma,\kg_{j})\;\,
         =\;\,\sum_{j=i+1}^m h^1(\Sigma,\kg_{j}).
       \end{array}       
       \as{\asf}
     \end{displaymath}

     \begin{varthm-roman}[Step 2]
       $h^1(\Delta_i,\kg_i)=h^0(\Delta_i,\kg_i)-
       \chi\Big(\ko_{\Delta_i}\big(D-\sum_{k=1}^{i-1}\Delta_k\big)\Big)
       +\deg(X_{i-1}\cap\Delta_i)$.
     \end{varthm-roman}
     We consider the exact sequence
     \begin{displaymath}
       \xymatrix@C=0.5cm{
         0\ar[r]&\kg_i
         \ar[r]&\ko_{\Delta_i}\left(D-\sum\nolimits_{k=1}^{i-1}\Delta_k\right)
         \ar[r]&\ko_{X_{i-1}\cap\Delta_i/\Delta_i}\left(D-\sum\nolimits_{k=1}^{i-1}\Delta_k\right)
         \ar[r]&0.
         }
     \end{displaymath}
     The result then follows from the long exact
     cohomology sequence.

     \begin{varthm-roman}[Step 3]
       $h^0\big(\Delta_i,\kg_i^0\big)-
       \chi\Big(\ko_{\Delta_i}\big(D-\sum_{k=1}^{i-1}\Delta_k\big)\Big)
       =h^1\big(\Delta_i,\kg_i^0\big)-\deg(X_i^0)$.
     \end{varthm-roman}
     This follows analogously, replacing $X_{i-1}$ by $X_i^0$, since
     $X_i^0=X_i^0\cap\Delta_i$.

     \begin{varthm-roman}[Step 4]
       $h^1\big(\Delta_i,\kg_i^0\big)\leq
       h^1\left(\Sigma,\kj_{X_i^0/\Sigma}\big(D-\sum_{k=1}^{i-1}\Delta_k\big)\right)=1$. 
     \end{varthm-roman}
     Note that $X_i^0:\Delta_i=\emptyset$, and hence
     $\kj_{X_i^0:\Delta_i/\Sigma}=\ko_\Sigma$. We thus have the following
     short exact sequence
     \begin{equation}\label{eq:irred-B:4}
       \xymatrix@C=0.6cm{
         0\ar[r]&
         \ko_\Sigma\left(D-\sum\nolimits_{k=1}^i\Delta_k\right)
         \ar[r]^(0.45){\cdot\Delta_i} & 
         \kj_{X_i^0/\Sigma}\left(D-\sum\nolimits_{k=1}^{i-1}\Delta_k\right)
         \ar[r] & 
         \kg_i^0\ar[r]&0.
         }
     \end{equation}
     By assumption (f) the divisor $D-K_\Sigma-\sum_{k=1}^i\Delta_k$
     is big and nef and hence
     \begin{displaymath}
       0=h^0\Big(\Sigma,\ko_\Sigma\big(-D+K_\Sigma+\mbox{$\sum_{k=1}^i$}\Delta_k\big)\Big)
       =h^2\Big(\Sigma,\ko_\Sigma\big(D-\mbox{$\sum_{k=1}^i$}\Delta_k\big)\Big).
     \end{displaymath}
     Thus the long exact cohomology sequence of \eqref{eq:irred-B:4}
     gives 
     \begin{displaymath}
       \xymatrix{
         H^1\Big(\Sigma,\kj_{X_i^0/\Sigma}\big(D-\sum_{k=1}^{i-1}\Delta_k\big)\Big)
         \ar[r]&
         H^1\big(\Delta_i,\kg_i^0\big)
         \ar[r]&
         0,
         }
     \end{displaymath}
     and 
     \begin{displaymath}
       h^1\big(\Delta_i,\kg_i^0\big)\leq
       h^1\left(\Sigma,\kj_{X_i^0/\Sigma}\big(D-\mbox{$\sum_{k=1}^{i-1}$}\Delta_k\big)\right).
     \end{displaymath}
     However, by assumption (b) the latter is just one.

     \begin{varthm-roman}[Step 5]
       $h^1(\Delta_i,\kg_i)\leq 1+\deg(X_{i-1}\cap\Delta_i)-\deg\big(X_i^0\big)$.
     \end{varthm-roman}
     We note that $\kg_i\hookrightarrow\kg_i^0$, and thus
       $h^0(\Delta_i,\kg_i)\leq h^0(\Delta_i,\kg_i^0\big)$.
     But then
     \begin{displaymath}
       \begin{array}{rcl}
         h^1(\Delta_i,\kg_i)&\leq_\expl{Step 2/3}&
         h^1\big(\Delta_i,\kg_i^0\big)-\deg\big(X_i^0\big)+\deg(X_{i-1}\cap\Delta_i)\\
         &\leq_\expl{Step 4}&
         1-\deg\big(X_i^0\big)+\deg(X_{i-1}\cap\Delta_i).
       \end{array}
     \end{displaymath}
     
     \begin{varthm-roman}[Step 6]
       Finish the proof.
     \end{varthm-roman}
     Taking into account, that
     $h^1(\Sigma,\kg_i)=h^1(\Delta_i,\kg_i)$, since $\kg_i$ is
     concentrated on $\Delta_i$\tom{ (cf.~\cite{Har77} III.2.10)},
     the first inequality follows from Step 1, while the second
     inequality is a consequence of Step 5 and the last inequality
     follows from assumption (c).
   \end{proof}

   In the Lemmata \ref{lem:irred-C}, \ref{lem:irred-D} and
   \ref{lem:irred-D*} we consider special classes of surfaces which
   allow us to do the necessary estimates in order to finally derive
   \begin{displaymath}
     \sum_{i=1}^m \big(\#
     X_i^0-\dim|\Delta_i|_l\big)>h^1\big(\Sigma,\kj_{X_0/\Sigma}(D)\big).      
   \end{displaymath}
   We first consider surfaces with Picard number one.

   \begin{lemma}\label{lem:irred-C}
     Let $\Sigma$ be a surface such that
     \begin{myenumerate}
     \item[\rm(i)] $\NS(\Sigma)=L\cdot\Z$ and $L$ ample, and
     \item[\rm(ii)] $h^1(\Sigma,C)=0$, whenever $C$ is effective.
     \end{myenumerate}

     Let $D\in\Div(\Sigma)$ and  $X_0\subset\Sigma$  a
     zero-dimensional scheme satisfying (0)--(3) from Lemma
     \ref{lem:irred-A} and
     \begin{myenumerate}
     \item[\rm(4)] $\sum\limits_{z\in\Sigma}\big(\deg(X_{0,z})\big)^2
       <
       \gamma\cdot (D-K_\Sigma)^2
       %=\gamma\cdot (d-\kappa)^2\cdot L^2
       $,\hspace*{0.5cm}
       where
       $\gamma=\tfrac{\big(1+\sqrt{1-4\beta}\big)^2\cdot
         L^2}{4\cdot\chi(\ko_\Sigma)+\max\{0,2\cdot K_\Sigma.L\}+6\cdot L^2}$.
     \end{myenumerate}

     Then, using the notation of Lemma \ref{lem:irred-A} and setting
     $X_S=\bigcup_{i=1}^m X_i^0$,
     \begin{displaymath}
       h^1\big(\Sigma,\kj_{X_0/\Sigma}(D)\big)+
       \sum_{i=1}^m
       \Big(h^0\big(\Sigma,\ko_\Sigma(\Delta_i)\big)-1\Big) < \# X_S.
     \end{displaymath}
   \end{lemma}
   \begin{proof}
     We fix the following notation: 
     \begin{displaymath}
       D\sim_a d\cdot L,\;\; K_\Sigma\sim_a\kappa\cdot
       L,\;\;\Delta_i\sim_a\delta_i\cdot L,\;\;\mbox{ and }
       l=\sqrt{L^2}>0.       
     \end{displaymath}
     Furthermore, we have
     $\gamma=\frac{\big(1+\sqrt{1-4\beta}\big)^2}{4\alpha}$, where 
     \begin{displaymath}
       \alpha
       =\tfrac{4\cdot \chi(\ko_\Sigma)+\max\{0,2\cdot K_\Sigma.L\}+6\cdot
         L^2}{4\cdot L^2}=
       \left\{
       \begin{array}[m]{ll}
         \frac{\chi(\ko_\Sigma)}{l^2}+\frac{\kappa+3}{2},&\mbox{ if }
         \kappa\geq 0,\\
         \frac{\chi(\ko_\Sigma)}{l^2}+\frac{3}{2},&\mbox{ if }
         \kappa< 0,
       \end{array}
       \right.
     \end{displaymath}

     \begin{varthm-roman}[Step 1]
       By (i) $\Sigma$ satisfies the assumption (*) of Lemma \ref{lem:irred-A}.
     \end{varthm-roman}
%     If $c\cdot L\sim_a C\subset\Sigma$ is effective, then in particular
%     $c=\frac{1}{l^2}\cdot C.L> 0$, and thus $C$ is ample, in particular
%     nef. Hence (*) in Lemma \ref{lem:irred-A} is fulfilled.

     \begin{varthm-roman}[Step 2]
       $\sum_{i=1}^m\delta_i\cdot l\leq \frac{(d-\kappa)\cdot l}{2}
       -\sqrt{\frac{(d-\kappa)^2\cdot l^2}{4}-\deg(X_S)}$, by \eqref{eq:irred-A:1}.
     \end{varthm-roman}

     \begin{varthm-roman}[Step 3]
       $h^1\big(\Sigma,\kj_{X_0}(D)\big)\leq
       \big(\kappa\cdot\sum_{i=1}^m\delta_i\big)\cdot l^2+
       \frac{1}{2}\Big(\big(\sum_{i=1}^m\delta_i\big)^2+\sum_{i=1}^m\delta_i^2\Big)\cdot l^2+m$.
     \end{varthm-roman}
     By Lemma \ref{lem:irred-B} we know:
     \begin{displaymath}
       \as{1.4}
       \begin{array}{rcl}
         h^1\big(\Sigma,\kj_{X_0}(D)\big)&\leq&
         \sum_{i=1}^m
         \Big(\Delta_i\cdot\big(K_\Sigma+\sum_{k=1}^i\Delta_k\big)+1\Big)\\
%         &=&
%         \sum_{i=1}^m
%         \delta_i\cdot\big(\kappa+\sum_{k=1}^i\delta_k\big)\cdot l^2+m\\
%         &=&
%         \big(\kappa\cdot \sum_{i=1}^m\delta_i\big)\cdot l^2+
%         \big(\sum_{i=1}^m\sum_{k=1}^i\delta_i\cdot\delta_k\big)\cdot
%         l^2+m\\
         &=&
         \big(\kappa\cdot\sum_{i=1}^m\delta_i\big)\cdot l^2+
         \frac{1}{2}\Big(\big(\sum_{i=1}^m\delta_i\big)^2+\sum_{i=1}^m\delta_i^2\Big)\cdot l^2+m. 
       \end{array}
       \as{\asf}
     \end{displaymath}

     \begin{varthm-roman}[Step 4] 
       $\sum_{i=1}^m
       \Big(h^0\big(\Sigma,\ko_\Sigma(\Delta_i)\big)-1\Big)
       \leq
       m\cdot\big(\chi(\ko_\Sigma)-1\big)+\frac{l^2}{2}\cdot\sum_{i=1}^m\delta_i^2
       -\frac{\kappa\cdot l^2}{2}\cdot\sum_{i=1}^m\delta_i$.
     \end{varthm-roman}
     Since $\Delta_i$ is effective by (ii), $h^1(\Sigma,\Delta_i)=0$. Hence by
     Riemann-Roch 
%     \begin{displaymath}
%       h^0(\Sigma,\Delta_i)\leq\chi\big(\ko_\Sigma(\Delta_i)\big)
%       =\frac{\Delta_i^2-K_\Sigma.\Delta_i}{2}+\chi(\ko_\Sigma).
%     \end{displaymath}
%     This implies
     \begin{displaymath}
       \as{1.4}
       \begin{array}{rcl}
         \sum_{i=1}^m
         \Big(h^0\big(\Sigma,\ko_\Sigma(\Delta_i)\big)-1\Big)
         &\leq&
         -m+m\cdot\chi(\ko_\Sigma)+\frac{1}{2}\sum_{i=1}^m\big(\Delta_i^2-K_\Sigma.\Delta_i\big)\\
         &=&m\cdot\big(\chi(\ko_\Sigma)-1\big)+\frac{l^2}{2}\cdot\sum_{i=1}^m\delta_i^2
         -\frac{\kappa\cdot l^2}{2}\cdot\sum_{i=1}^m\delta_i.
       \end{array}
       \as{\asf}
     \end{displaymath}

     \begin{varthm-roman}[Step 5]       
       Finish the proof.
     \end{varthm-roman}
     In the following consideration we use that $\deg(X_S)\leq
     \deg(X_0)\leq\beta\cdot (d-\kappa)^2\cdot l^2$. 
     \begin{displaymath}
       \as{2}
       \begin{array}{l}
         h^1\big(\Sigma,\kj_{X_0}(D)\big)+
         \sum_{i=1}^m
         \Big(h^0\big(\Sigma,\ko_\Sigma(\Delta_i)\big)-1\Big)\\
         \leq_\expl{Step 3 / 4}\;\,
         m\cdot\chi(\ko_\Sigma)+
         l^2\cdot\sum_{i=1}^m\delta_i^2
         +\frac{\kappa\cdot l^2}{2}\cdot\sum_{i=1}^m\delta_i
         +\frac{l^2}{2}\cdot\big(\sum_{i=1}^m\delta_i\big)^2\\
         \leq\;\,
         \alpha\cdot\big(l\cdot\sum_{i=1}^m\delta_i\big)^2\;\,
         \leq_\expl{Step 2}\;\,
         \alpha\cdot
         \left(\frac{(d-\kappa)\cdot l}{2}
           -\sqrt{\frac{(d-\kappa)^2\cdot l^2}{4}-\deg(X_S)}\right)^2\\
         \leq\;\,\alpha\cdot
         \left(\frac{\frac{(d-\kappa)^2\cdot l^2}{4}
             -\big(\frac{(d-\kappa)^2\cdot l^2}{4}-\deg(X_S)\big)}
           {\frac{(d-\kappa)\cdot l}{2}
             +\sqrt{\frac{(d-\kappa)^2\cdot l^2}{4}-\deg(X_S)}}
         \right)^2\;\,
         =\;\,
         \alpha\cdot
         \left(\frac{2\cdot\deg(X_S)\big)}
           {(d-\kappa)\cdot l
             +\sqrt{(d-\kappa)^2\cdot l^2-4\cdot\deg(X_S)}}
         \right)^2\\
         \leq\;\,
         \frac{4\alpha}{\big(1+\sqrt{1-4\beta}\big)^2\cdot(d-\kappa)^2\cdot l^2}\cdot
         \big(\deg(X_S)\big)^2\;\,
         =\;\,
         \frac{1}{\gamma\cdot (D-K_\Sigma)^2}\cdot
         \big(\sum_{z\in\Sigma}\deg(X_{S,z})\big)^2\\
         \leq\;\,
         \frac{\# X_S}{\gamma\cdot (D-K_\Sigma)^2}\cdot
         \sum_{z\in\Sigma}\deg(X_{S,z})^2\;\,
         \leq\;\,
         \frac{\# X_S}{\gamma\cdot (D-K_\Sigma)^2}\cdot
         \sum_{z\in\Sigma}\deg(X_{0,z})^2\;\,
         <_\expl{(4)}\;\,\# X_S.
       \end{array}
       \as{\asf}
     \end{displaymath}
   \end{proof}

   The second class of surfaces which we consider, are products of curves.
%   We use the notation of Section
%   \ref{subsec:product-curves}. 

   \begin{lemma}\label{lem:irred-D}
     Let $C_1$ and $C_2$ be two smooth projective curves of genera $g_1$
     and $g_2$  respectively with $g_1\geq g_2\geq 0$, such that for $\Sigma=C_1\times C_2$
     the N\'eron--Severi group is $\NS(\Sigma)=C_1\Z\oplus C_2\Z$,
     and let $D\in\Div(\Sigma)$ such that $D\sim_a aC_1+bC_2$ with
     $a>\max\{2g_2-2,2-2g_2\}$ and $b>\max\{2g_1-2,2-2g_1\}$. 
     Suppose moreover that $X_0\subset\Sigma$ is a
     zero-dimensional scheme satisfying (1)--(3) from Lemma
     \ref{lem:irred-A} and
     \begin{enumerate}
%     \item[(1)] $\exists\;C_0\in|D|_l\;\mbox{ irreducible}:\;X_0\subset C_0$,
%     \item[(2)] $h^1\big(\Sigma,\kj_{X_0/\Sigma}(D)\big)>0$, and 
%     \item[(3)] $\deg(X_0)<\beta\cdot (D-K_\Sigma)^2$ for some $0<\beta\leq
%       \frac{1}{4}$.
     \item[(4)] $\sum\limits_{z\in\Sigma}\big(\deg(X_{0,z})\big)^2
       \;<\;\gamma\cdot(D-K_\Sigma)^2$, 
     \end{enumerate}
     where $\gamma$ may be taken from the table in Theorem \ref{thm:irred-products-of-curves}.
%     with $\alpha=\frac{a-2g_2+2}{b-2g_1+2}>0$.
%     \begin{center}
%       \tablefirsthead{\hline
%         \multicolumn{1}{|c}{$g_1$}&\multicolumn{1}{|c|}{$g_2$}&$\gamma$\\\hline\hline}
%       \tablehead{\hline
%         \multicolumn{1}{|c}{$g_1$}&\multicolumn{1}{|c|}{$g_2$}&$\gamma$\\\hline\hline}
%       \tabletail{\hline}
%       \tablelasttail{\\\hline}
%       \begin{supertabular}{|r|r|c|}         
%         $0$&$0$&$\frac{1}{24} $\\
%         $1$&$0$&$\frac{1}{\max\{32,2\alpha\}}$\\
%         $\geq 2$&$0$&
%         $\frac{1}{\max\{24+16g_1,4g_1\alpha\}} $\\
%         $1$&$1$&
%         $\frac{1}{\max\big\{32,2\alpha,\tfrac{2}{\alpha}\big\}} $\\
%         $\geq 2$&$\geq 1$&
%         $\frac{1}{\max\left\{24+16g_1+16g_2,4g_1\alpha,\tfrac{4g_2}{\alpha}\right\}}$
%       \end{supertabular}
%     \end{center}

     Then, using the notation of Lemma \ref{lem:irred-A} and
     setting\tom{\footnote{Remember that $\# X_S$ is the number of  
         points in the support of $X_S$.}}
     $X_S=\bigcup_{i=1}^m X_i^0$, 
     \begin{displaymath}
       h^1\big(\Sigma,\kj_{X_0}(D)\big)+
       \sum_{i=1}^m
       \Big(h^0\big(\Sigma,\ko_\Sigma(\Delta_i)\big)-1\Big) < \# X_S.
     \end{displaymath}
   \end{lemma}
   \begin{proof}
     Then $K_\Sigma\sim_a(2g_2-2)\cdot C_1+(2g_1-2)\cdot C_2$ and we fix the notation:
     \begin{displaymath}
       \Delta_i\sim_a a_iC_1+b_iC_2,\;\;\; \kappa_1=a-2g_2+2\;\;\; \text{
         and }\;\;\; \kappa_2=b-2g_1+2. 
     \end{displaymath}

     \begin{varthm-roman}[Step 1]
       $\Sigma$ satisfies the assumption (*) of Lemma
       \ref{lem:irred-A}. Moreover, due to the
       assumptions on $a$ and $b$ we know that $D-K_\Sigma$ is ample
       and $D+K_\Sigma$ is nef, i.~e.~(0) in Lemma \ref{lem:irred-A}
       is fulfilled as well.
     \end{varthm-roman}

     \begin{varthm-roman}[Step 2a]
       $\big(\frac{\kappa_1}{4}\big)\cdot\sum_{i=1}^m
       b_i+\big(\frac{\kappa_2}{4}\big)\cdot\sum_{i=1}^m a_i 
       \leq \deg(X_S)$.
     \end{varthm-roman}
     Let us first notice that the strict inequality ``$<$'' in Lemma
     \ref{lem:irred-A} (e) for ample divisors $H$ comes down to ``$\leq$''
     for nef divisors $H$. We may apply this for $H=C_1$ and $H=C_2$ and
     deduce the following inequalities:
     \begin{equation}
       \label{eq:irred-D:1}
       0\leq \left(D-K_\Sigma-\sum_{k=1}^i\Delta_k-\Delta_i\right).C_1
       =\kappa_2-\sum_{k=1}^i b_k-b_i,
     \end{equation}
     and
     \begin{equation}
       \label{eq:irred-D:2}
       0\leq \left(D-K_\Sigma-\sum_{k=1}^i\Delta_k-\Delta_i\right).C_2
       =\kappa_1-\sum_{k=1}^i a_k-a_i.
     \end{equation}
     For the following consideration we choose
     $i_0,j_0\in\{1,\ldots,m\}$ such that $a_{i_0}\geq a_i$ for all
     $i=1,\ldots,m$ and $b_{j_0}\geq b_j$ for all
     $j=1,\ldots,m$. Then
     \begin{equation}
       \label{eq:irred-D:3}
       \kappa_1\geq 2a_i\;\;\; \text{ and }\;\;\; \kappa_2\geq 2b_j
%       \kappa_1\geq\sum_{k=1}^{i_0}a_k+a_{i_0}\geq 2a_{i_0}\geq 2a_i 
%       \intertext{ and }
%       \label{eq:irred-D:3*}
%       \kappa_2\geq\sum_{k=1}^{j_0}b_k+b_{j_0}\geq 2b_{j_0}\geq 2b_j 
     \end{equation}
     for all $i,j=1,\ldots,m$;  finally
     \eqref{eq:irred-D:1}--\eqref{eq:irred-D:3} lead to
     \begin{displaymath}
       \as{1.6}
       \begin{array}{l}
         \deg(X_S)\;\,=\;\,\sum_{i=1}^m\deg\big(X_i^0\big)
         \;\,\geq_\expl{Lemma \ref{lem:irred-A} (c)}
         \;\,\sum_{i=1}^m\big(D-K_\Sigma-\sum_{k=1}^i\Delta_k\big).\Delta_i\\
%         \;\,=\;\,
%         (D-K_\Sigma).\sum_{i=1}^m\Delta_i-\sum_{1\leq k\leq i\leq
%           1}\Delta_k.\Delta_i\\
%         \;\,=\;\,(D-K_\Sigma).\sum_{i=1}^m\Delta_i
%         -\frac{1}{2}\sum_{i=1}^m\Delta_i^2-\frac{1}{2}\big(\sum_{i=1}^m\Delta_i\big)^2\\
         \;\,=\;\,
         \kappa_1 \sum_{i=1}^m b_i +\kappa_2 \sum_{i=1}^m a_i
         -\sum_{i=1}^m a_ib_i
         %\\\hspace*{1cm}
         -\sum_{i=1}^ma_i \sum_{i=1}^mb_i\\
         \;\,\geq\;\,%_\expl{\eqref{eq:irred-D:1} / \eqref{eq:irred-D:2}}\;\,
%         \frac{\kappa_1}{2} \sum_{i=1}^m b_i+\frac{\kappa_2}{2} \sum_{i=1}^m a_i+
%         \frac{1}{2}\big(\sum_{i=1}^ma_i+a_m\big) \cdot\sum_{i=1}^m b_i\\
%         \;\;\;\;\;\;\;\;\;\;+\frac{1}{2}\big(\sum_{i=1}^mb_i+b_m\big) \cdot\sum_{i=1}^m a_i
%         -\sum_{i=1}^m
%         a_ib_i-\sum_{i=1}^ma_i \sum_{i=1}^mb_i\\
%         \;\,=\;\,
         \frac{\kappa_1}{2} \sum_{i=1}^m b_i+\frac{\kappa_2}{2} \sum_{i=1}^m a_i
         +\frac{a_m}{2} \sum_{i=1}^m b_i+\frac{b_m}{2} \sum_{i=1}^m a_i
         -\sum_{i=1}^m a_ib_i\\
         \;\,\geq\;\,%_\expl{\eqref{eq:irred-D:3} / \eqref{eq:irred-D:3*}}\;\,
%         \frac{\kappa_1}{4} \sum_{i=1}^m b_i+\frac{1}{4}\sum_{i=1}^m 2a_ib_i
%         +\frac{\kappa_2}{4} \sum_{i=1}^m a_i+\frac{1}{4}\sum_{i=1}^m 2a_ib_i\\
%         \;\;\;\;\;\;\;\;\;\;\;\;\;\;+\frac{a_m}{2} \sum_{i=1}^m b_i+\frac{b_m}{2} \sum_{i=1}^m a_i
%         -\sum_{i=1}^m a_ib_i\\
%         \;\,\geq\;\,
         \frac{\kappa_1}{4} \sum_{i=1}^m
         b_i+\frac{\kappa_2}{4} \sum_{i=1}^m a_i.
       \end{array}
       \as{\asf}
     \end{displaymath}

     \begin{varthm-roman}[Step 2b]
       $\sum_{i=1}^ma_i\cdot\sum_{i=1}^mb_i
       \leq\frac{8}{(D-K_\Sigma)^2}\cdot\big(\deg(X_S)\big)^2$.
     \end{varthm-roman}
     Using Step 2a we deduce
     \begin{displaymath}
       \as{1.6}
       \begin{array}{rcl}
         \big(\deg(X_S)\big)^2 &>&
         \Big(\frac{\kappa_2}{4}\cdot\sum_{i=1}^m a_i +\frac{\kappa_1}{4}\cdot\sum_{i=1}^m b_i\Big)^2\\
         &\geq&
         \frac{4\cdot\kappa_1\cdot\kappa_2}{16}\cdot\sum_{i=1}^m a_i\cdot\sum_{i=1}^m b_i\\
         &=&
         \frac{(D-K_\Sigma)^2}{8}\cdot\sum_{i=1}^m a_i\cdot\sum_{i=1}^m b_i.
       \end{array}
       \as{\asf}
     \end{displaymath}

     \begin{varthm-roman}[Step 2c]
       $\sum_{i=1}^ma_i\leq
       \left\{
         \begin{array}{ll}
           \frac{2\alpha}{(D-K_\Sigma)^2}\cdot\big(\deg(X_S)\big)^2,&
           \text{ if } \sum_{i=1}^m b_i=0,\\
           \frac{8}{(D-K_\Sigma)^2}\cdot\big(\deg(X_S)\big)^2,&
           \text{ otherwise}.
         \end{array}
       \right.$ 
     \end{varthm-roman} 
     If $\sum_{i=1}^m b_i=0$, then the same consideration as in Step 2a shows
     \begin{displaymath}
       \deg(X_S)\geq \kappa_2\cdot\sum_{i=1}^ma_i>0,
     \end{displaymath}
     and thus
     \begin{displaymath}
       \tfrac{(D-K_\Sigma)^2}{2\alpha}\cdot\sum_{i=1}^ma_i
       \leq
       \kappa_2^2\cdot\left(\sum_{i=1}^ma_i\right)^2\leq\big(\deg(X_S)\big)^2.
     \end{displaymath}
     If $\sum_{i=1}^m b_i\not=0$, then we are done by Step 2b.
%     \begin{displaymath}
%       \sum_{i=1}^ma_i\leq
%       \sum_{i=1}^ma_i\cdot\sum_{i=1}^mb_i
%       \leq
%       \tfrac{8}{(D-K_\Sigma)^2}\cdot\big(\deg(X_S)\big)^2.
%     \end{displaymath}

     \begin{varthm-roman}[Step 2d]
       $\sum_{i=1}^mb_i\leq
       \left\{
         \begin{array}{ll}
           \frac{2}{\alpha\cdot(D-K_\Sigma)^2}\cdot\big(\deg(X_S)\big)^2,&
           \text{ if } \sum_{i=1}^m a_i=0,\\
           \frac{8}{(D-K_\Sigma)^2}\cdot\big(\deg(X_S)\big)^2,&
           \text{ otherwise.}
         \end{array}
       \right.$ 
     \end{varthm-roman} 
     This is proved in the same way as Step 2c.

     \begin{varthm-roman}[Step 3]
       $h^1\big(\Sigma,\kj_{X_0}(D)\big)\leq
       2\sum\limits_{i=1}^ma_i\sum\limits_{i=1}^mb_i
       +(2g_1-2)\sum\limits_{i=1}^ma_i+(2g_2-2)\sum\limits_{i=1}^mb_i+m$.
     \end{varthm-roman}
     The following sequence of inequalities is due to 
     Lemma \ref{lem:irred-B} and the fact that $\Delta_i.\Delta_j\geq 0$
     for any $i,j\in\{1,\ldots,m\}$:
     \begin{displaymath}
       \as{1.4}
       \begin{array}{l}
         h^1\big(\Sigma,\kj_{X_0}(D)\big)\;\,\leq\;\,
         \sum_{i=1}^m
         \Big(\Delta_i\cdot\big(K_\Sigma+\sum_{k=1}^i\Delta_k\big)+1\Big)\\
%         \;\,=\;\,
%         K_\Sigma\cdot \sum_{i=1}^m\Delta_i +
%         \sum_{1\leq k\leq i\leq m}\Delta_i.\Delta_k+m\\
         \;\,\leq\;\,
         K_\Sigma\cdot \sum_{i=1}^m\Delta_i +
         \big(\sum_{i=1}^m\Delta_i\big)^2+m \\
         \;\,=\;\,
         (2g_1-2)\cdot \sum_{i=1}^ma_i+(2g_2-2)\cdot \sum_{i=1}^mb_i
         +2\cdot \sum_{i=1}^ma_i\cdot\sum_{i=1}^mb_i+m.\\
       \end{array}
       \as{\asf}
     \end{displaymath}

     \begin{varthm-roman}[Step 4]
       We find the  estimate $\sum_{i=1}^m\Big(h^0\big(\Sigma,\ko_\Sigma(\Delta_i)\big)-1\Big)
       \leq\beta$, where
       \begin{displaymath}
         \beta
         =
         \left\{
           \begin{array}{l}
             \sum_{i=1}^ma_i\cdot\sum_{i=1}^mb_i
             +\sum_{i=1}^mb_i,\;\;\;\text{ if } g_1= 1,g_2=0,\\
             \sum_{i=1}^ma_i\cdot\sum_{i=1}^mb_i-m,\;\;\;\text{ if }
             g_1=1,g_2=1,\exists\;i_0\;:\;a_{i_0}b_{i_0}>0,\\
             \sum_{i=1}^ma_i+\sum_{i=1}^mb_i-m,\;\;\;\text{ if }
             g_1=1,g_2=1,\forall\;i\;:\;a_ib_i=0,\\
             \sum_{i=1}^ma_i\cdot\sum_{i=1}^mb_i
             +\sum_{i=1}^ma_i+\sum_{i=1}^mb_i,\;\;\; \text{ otherwise.} \\
           \end{array}
         \right.
       \end{displaymath}
     \end{varthm-roman}

     In general \tom{by
     Corollary \ref{cor:kuenneth-formula}}
     $h^0\big(\Sigma,\ko_\Sigma(\Delta_i)\big)\leq a_ib_i+a_i+b_i+1$, whereas
     if  $g_1=1,g_2=0$ \tom{by Lemma \ref{lem:divsors-on-ruled-surfaces}} we have
     $h^0\big(\Sigma,\ko_\Sigma(\Delta_i)\big)= a_ib_i+b_i+1$. It thus only
     remains to consider the case $g_1=g_2=1$, where we get\tom{\footnote{Applying Lemma
       \ref{lem:kuenneth-formula-elliptic-curves}.}}     
     \begin{displaymath}
       \sum_{i=1}^mh^0\big(\Sigma,\ko_\Sigma(\Delta_i)\big) 
       =
       \sum_{a_i,b_i>0}a_ib_i+\sum_{a_i=0} b_i+\sum_{b_i=0} a_i.   
     \end{displaymath}
     If always either $a_i$ or $b_i$ is zero, we are done.
     Otherwise there exists some $i_0\in\{1,\ldots,m\}$ such that
     $a_{i_0}\not=0\not=b_{i_0}$. Then looking at the right hand side
     we see
     \begin{displaymath}
       \sum_{i=1}^mh^0\big(\Sigma,\ko_\Sigma(\Delta_i)\big) 
       \leq
       \sum_{a_i,b_i>0}a_ib_i+a_{i_0}\cdot\sum_{a_i=0} b_i
       +b_{i_0}\cdot\sum_{b_i=0} a_i
       \leq \sum_{i=1}^ma_i\cdot\sum_{i=1}^mb_i.
     \end{displaymath}

     \begin{varthm-roman}[Step 5]       
       Finish the proof.
     \end{varthm-roman}

     Using Step 3 and Step 4, and taking $m\leq\sum_{i=1}^ma_i+b_i$
     into account, we get $h^1\big(\Sigma,\kj_{X_0}(D)\big)+
     \sum_{i=1}^m
     \Big(h^0\big(\Sigma,\ko_\Sigma(\Delta_i)\big)-1\Big)\leq
     \beta'$, where $\beta'$ may be chosen as
     \begin{displaymath}
       \beta'
       =
       \left\{
         \begin{array}{ll}
           3\cdot\sum_{i=1}^ma_i\cdot\sum_{i=1}^mb_i,
           &\text{if } g_1= 0,g_2=0,\\
           3\cdot \sum_{i=1}^ma_i\cdot\sum_{i=1}^mb_i
           +\sum_{i=1}^ma_i,&\text{if } g_1= 1,g_2=0,\\
%           3\cdot\sum_{i=1}^ma_i\cdot\sum_{i=1}^mb_i
%           +2g_1\cdot\sum_{i=1}^ma_i,&\text{if } g_1\geq 2,g_2=0, \\
           3\cdot\sum_{i=1}^ma_i\cdot\sum_{i=1}^mb_i
           +2g_1\cdot\sum_{i=1}^ma_i+2g_2\cdot\sum_{i=1}^mb_i,&\text{if } g_1\geq 2,g_2\geq 0.
         \end{array}
       \right.
     \end{displaymath}
     For the case $g_1=g_2=1$ we take a closer look. We find at once the following upper bounds
     $\beta''$ for $h^1\big(\Sigma,\kj_{X_0}(D)\big)+
     \sum_{i=1}^m
     \Big(h^0\big(\Sigma,\ko_\Sigma(\Delta_i)\big)-1\Big)$
     \begin{displaymath}
       \beta''
       =
       \left\{
         \begin{array}{ll}
           3\cdot \sum_{i=1}^ma_i\cdot\sum_{i=1}^mb_i,
           &\;\text{if } \exists\;i_0\;:\;a_{i_0}b_{i_0}\not=0,\\
           2\cdot\sum_{i=1}^ma_i\cdot\sum_{i=1}^mb_i+\sum_{i=1}^ma_i+\sum_{i=1}^mb_i,
           &\;\text{if } \forall\;i:a_ib_i=0.
         \end{array}
       \right.       
     \end{displaymath}
     Considering now the cases $\sum_{i=1}^ma_i\not=0
     \not=\sum_{i=1}^mb_i$, $\sum_{i=1}^ma_i=0$ and
     $\sum_{i=1}^mb_i=0$, we can replace these by 
     \begin{displaymath}
       \beta''\leq\beta'
       =
       \left\{
         \begin{array}{ll}
           4\cdot \sum_{i=1}^ma_i\cdot\sum_{i=1}^mb_i,
           &\text{if } \sum_{i=1}^ma_i\not=0\not=\sum_{i=1}^mb_i,\\
           \sum_{i=1}^ma_i,
           &\text{if } \sum_{i=1}^mb_i=0,\\
           \sum_{i=1}^mb_i,
           &\text{if } \sum_{i=1}^ma_i=0.
         \end{array}
       \right.       
     \end{displaymath}

     Applying now the results of Step 2 in all cases we get
     \begin{displaymath}
       \as{1.6}
       \begin{array}{l}
         h^1\big(\Sigma,\kj_{X_0}(D)\big)+
         \sum_{i=1}^m
         \Big(h^0\big(\Sigma,\ko_\Sigma(\Delta_i)\big)-1\Big)
         \;\,\leq\;\,\beta'\;\,\leq\;\,
         \frac{1}{\gamma\cdot(D-K_\Sigma)^2}\cdot\big(\deg(X_S)\big)^2\\
         \hspace*{1cm}=\;\,
         \frac{1}{\gamma\cdot(D-K_\Sigma)^2}\cdot
         \big(\sum_{z\in\Sigma}\deg(X_{S,z})\big)^2
         \;\,\leq\;\,
         \frac{\# X_S}{\gamma\cdot(D-K_\Sigma)^2}\cdot
         \sum_{z\in\Sigma}\deg(X_{S,z})^2\\
         \hspace*{1cm}\leq\;\,
         \frac{\# X_S}{\gamma\cdot(D-K_\Sigma)^2}\cdot
         \sum_{z\in\Sigma}\deg(X_{0,z})^2\;\,
         <_\expl{(4)}\;\,\# X_S.
       \end{array}
       \as{\asf}
     \end{displaymath}
   \end{proof}

   \begin{remark}
     Lemma \ref{lem:irred-D}, and hence Theorem
     \ref{thm:irred-products-of-curves} could easily be generalised to
     other surfaces $\Sigma$ with irreducible curves $C_1,C_2\subset
     \Sigma$ such that 
       $\NS(\Sigma)=C_1\Z\oplus C_2\Z$ with intersection matrix 
       $\left(
       \begin{smallmatrix}
         0&1\\1&0
       \end{smallmatrix}
       \right)$
     once we have an estimate similar to 
     \begin{displaymath}
       h^0(\Sigma, aC_1+bC_2)\leq ab+a+b+1
     \end{displaymath}
     for an effective divisor $aC_1+bC_2$. 

     With a number of small modifications we are even able to adapt it
     in the following lemma
     in the case of geometrically ruled surfaces with non-positive
     invariant $e$ although the intersection pairing looks more
     complicated. 

     The problem with arbitrary geometrically ruled surfaces is
     the existence of the section with negative
     self-intersection, once the invariant $e>0$, since then the proof
     of Lemma \ref{lem:irred-A} no longer works.
   \end{remark}

   In the following lemma we use the notation of Subsection
   \ref{subsec:irreducibility:ruled-surfaces}.

   \begin{lemma}\label{lem:irred-D*}
     Let $\pi:\Sigma\rightarrow C$ be a geometrically ruled surface
     with invariant $e\leq 0$ and $g=g(C)$, 
     and let $D\in\Div(\Sigma)$ such that $D\sim_a aC_0+bF$ with $a\geq 2$, $b>2g-2+\frac{ae}{2}$, and
     if $g=0$ then $b\geq 2$. 
     Suppose moreover that $X_0\subset\Sigma$ is a
     zero-dimensional scheme satisfying (1)--(3) from Lemma
     \ref{lem:irred-A} and
     \begin{enumerate}
%     \item[(1)] $\exists\;C_0\in|D|_l\;\mbox{ irreducible}:\;X_0\subset C_0$,
%     \item[(2)] $h^1\big(\Sigma,\kj_{X_0/\Sigma}(D)\big)>0$, and 
%     \item[(3)] $\deg(X_0)<\beta\cdot (D-K_\Sigma)^2$ for some $0<\beta\leq
%       \frac{1}{4}$.
     \item[(4)] $\sum\limits_{z\in\Sigma}\big(\deg(X_{0,z})\big)^2
       \;<\;\gamma\cdot(D-K_\Sigma)^2$, 
     \end{enumerate}
     where $\gamma$ may be taken from the  table in Theorem \ref{thm:irred-ruled-surfaces}.
%     with
%     $\alpha=\frac{a+2}{b+2-2g-\tfrac{ae}{2}}>0$. 
%     \begin{displaymath}
%       \begin{array}{|r|r|c|}
%         \hline
%         \multicolumn{1}{|c}{g}&\multicolumn{1}{|c|}{e}&\gamma\\\hline\hline
%%         g=0\text{ or }g=1& \frac{1}{24}&\frac{1}{48} \\
%%         g\geq 2&\frac{1}{16+16g} &
%%         \frac{1}{48+(2g-1)\cdot\max\left\{2\alpha,16\right\}}\\\hline 
%          0 & 0 &  \tfrac{1}{24}\\
%          1 & 0 &  \tfrac{1}{\max\{24,2\alpha\}}\\
%          1 & -1 & \tfrac{1}{\max\Big\{\min\big\{30+\tfrac{16}{\alpha}+4\alpha,40+9\alpha\big\},
%           \tfrac{13}{2}\alpha\Big\}}\\
%          \geq 2 & 0 &  \tfrac{1}{\max\{24+16g, 4g\alpha\}}\\
%          \geq 2 & <0 & 
%           \tfrac{1}{\max\Big\{\min\big\{24+16g-9e\alpha,18+16g-9e\alpha
%           -\tfrac{16}{e\alpha}\big\},4g\alpha-9e\alpha\Big\}}\\\hline
%       \end{array}
%     \end{displaymath}

     Then, using the notation of Lemma \ref{lem:irred-A} and
     setting\tom{\footnote{Remember that $\# X_S$ is the number of  
         points in the support of $X_S$.}} 
     $X_S=\bigcup_{i=1}^m X_i^0$, 
     \begin{displaymath}
       h^1\big(\Sigma,\kj_{X_0}(D)\big)+
       \sum_{i=1}^m
       \Big(h^0\big(\Sigma,\ko_\Sigma(\Delta_i)\big)-1\Big) < \# X_S.
     \end{displaymath}
   \end{lemma}
   \begin{proof}
     Remember that the N\'eron--Severi group of $\Sigma$ is generated
     by a section $C_0$ of $\pi$ and a fibre $F$ with intersection
     pairing given by
     $\left(\begin{smallmatrix}
         -e&1\\1&0 
       \end{smallmatrix}\right)$. 
     Then $K_\Sigma\sim_a-2C_0+(2g-2-e)\cdot F$ and we fix the notation:
     \begin{displaymath}
       \Delta_i\sim_a a_iC_0+b_i'F.
     \end{displaymath}
     Note that then
     \begin{displaymath}
       a_i\geq 0\;\;\;\; \text{ and }\;\;\;\; b_i:=b_i'-\tfrac{e}{2}a_i\geq 0.
     \end{displaymath}
     Finally we set $\kappa_1=a+2$ and
     $\kappa_2=b+2-2g-\tfrac{ae}{2}$ and get
     \begin{equation}\label{eq:irred-D*:0}
       (D-K_\Sigma)^2=-e\cdot(a+2)^2+2\cdot(a+2)\cdot(b+2+e-2g)=2\cdot\kappa_1\cdot\kappa_2.
     \end{equation}

     Replacing the equations \eqref{eq:irred-D:1} and
     \eqref{eq:irred-D:2} by 
     \begin{equation}
       \label{eq:irred-D*:1}
       0\leq \left(D-K_\Sigma-\sum_{k=1}^i\Delta_k-\Delta_i\right).\left(C_0+\tfrac{e}{2}F\right)
       =\kappa_2-\sum_{k=1}^i b_k-b_i,
     \end{equation}
     and
     \begin{equation}
       \label{eq:irred-D*:2}
       0\leq \left(D-K_\Sigma-\sum_{k=1}^i\Delta_k-\Delta_i\right).F
       =\kappa_1-\sum_{k=1}^i a_k-a_i,
     \end{equation}
     the assertions of Step 1 to Step 2c in the proof of Lemma
     \ref{lem:irred-D} remain literally true.

     \begin{varthm-roman}[Step 2d]
       $\left(\sum\limits_{i=1}^ma_i\right)^2
       \leq\frac{32\alpha}{(D-K_\Sigma)^2}\big(\deg(X_S)\big)^2$
       and
       $\left(\sum\limits_{i=1}^mb_i\right)^2
       \leq\frac{32}{\alpha\cdot(D-K_\Sigma)^2}\big(\deg(X_S)\big)^2$.
     \end{varthm-roman}
     This follows from the following inequality with the aid of Step 2a and \eqref{eq:irred-D*:0}, 
     \begin{multline*}
       \big(\deg(X_S)\big)^2 \geq
       \big(\tfrac{\kappa_2}{4}\cdot\mbox{$\sum\nolimits_{i=1}^m$} a_i\big)^2+
       \big(\tfrac{\kappa_1}{4}\cdot\mbox{$\sum\nolimits_{i=1}^m$} b_i\big)^2\\
       \geq
       \tfrac{2\cdot\kappa_1\cdot\kappa_2}{32\alpha}\cdot\big(\mbox{$\sum\nolimits_{i=1}^m$} a_i\big)^2+
       \tfrac{2\cdot\kappa_1\cdot\kappa_2\cdot\alpha}{32}\cdot\big(\mbox{$\sum\nolimits_{i=1}^m$} b_i\big)^2.
     \end{multline*}

     \begin{varthm-roman}[Step 3]
       $h^1\big(\Sigma,\kj_{X_0}(D)\big)\leq
       2\cdot\sum\limits_{i=1}^ma_i\cdot\sum\limits_{i=1}^mb_i
       +(2g-2)\cdot\sum\limits_{i=1}^ma_i-2\cdot\sum\limits_{i=1}^mb_i+m$
       is proved as Step 3 in Lemma \ref{lem:irred-D}.
     \end{varthm-roman}
     
%     By Lemma \ref{lem:irred-B} and since $\Delta_i.\Delta_j\geq 0$
%     for any $i,j\in\{1,\ldots,m\}$ we have:
%     \begin{displaymath}
%       \as{1.4}
%       \begin{array}{l}
%         h^1\big(\Sigma,\kj_{X_0}(D)\big)\;\,\leq\;\,
%         \sum_{i=1}^m
%         \Big(\Delta_i\cdot\big(K_\Sigma+\sum_{k=1}^i\Delta_k\big)+1\Big)\\
%         \;\,\leq\;\,
%%         K_\Sigma\cdot \sum_{i=1}^m\Delta_i +
%%         \big(\sum_{i=1}^m\Delta_i\big)^2 +m\\
%%         \;\,=\;\,
%         (2g-2)\cdot \sum_{i=1}^ma_i-2\cdot \sum_{i=1}^mb_i
%         +2\cdot \sum_{i=1}^ma_i\cdot\sum_{i=1}^mb_i+m.\\
%       \end{array}
%       \as{\asf}
%     \end{displaymath}

     \begin{varthm-roman}[Step 4a]
       If $e=0$, we find the  estimate 
       \begin{displaymath}
         \sum_{i=1}^m\Big(h^0\big(\Sigma,\ko_\Sigma(\Delta_i)\big)-1\Big)%\beta
         \leq
         \as{1.8}
         \left\{
           \begin{array}{ll}
             \sum\limits_{i=1}^ma_i\cdot\sum\limits_{i=1}^mb_i+\sum\limits_{i=1}^mb_i-m,&\text{if } g=1, 
             \sum_{i=1}^mb_i\not=0,\\
             \sum\limits_{i=1}^ma_i\cdot\sum\limits_{i=1}^mb_i+\sum\limits_{i=1}^mb_i=0,&\text{if } g=1, 
             \sum_{i=1}^mb_i=0,\\
             \sum\limits_{i=1}^ma_i\cdot\sum\limits_{i=1}^mb_i
             +\sum\limits_{i=1}^ma_i+\sum\limits_{i=1}^mb_i,&\text{for } g \text{
               arbitrary.} \\
           \end{array}
         \right.
         \as{\asf}
       \end{displaymath}
     \end{varthm-roman} 
     We note that in this case $b_i'=b_i$ and that $b_i=0$ thus implies $a_i>0$.
     But then
     \begin{displaymath}
       h^0\big(\Sigma,\ko_\Sigma(\Delta_i)\big)\leq 
       \left\{
         \begin{array}{ll}
           a_ib_i+b_i,&
           \text{ if } g=1, b_i>0,\\     
           a_ib_i+b_i+1=1,&
           \text{ if } g=1, b_i=0,\\     
           a_ib_i+a_i+b_i+1,&
           \text{ otherwise.}
         \end{array}
       \right.
       \end{displaymath}
       The results for $g$ arbitrary respectively  $g=1$ and
       $\sum_{i=1}^mb_i=0$ thus follow right away.  
       If, however, some $b_{i_0}>0$, then 
         $\sum_{i\not= j}a_ib_j\geq b_{i_0}\sum_{i\not=i_0}a_i\geq \#\big\{b_i\;|\;b_i=0\big\}$
%       \begin{displaymath}
%         \sum_{i\not= j}a_ib_j\geq b_{i_0}\sum_{i\not=i_0}a_i\geq \#\big\{b_i\;|\;b_i=0\big\}
%       \end{displaymath}
       and hence
       \begin{multline*}
         h^0\big(\Sigma,\ko_\Sigma(\Delta_i)\big)\leq 
         \sum_{i=1}^m a_ib_i+\sum_{i=1}^m b_i
         +\#\big\{b_i\;|\;b_i=0\big\}\\
         =
         \sum_{i=1}^m a_i\cdot\sum_{i=1}^mb_i+\sum_{i=1}^m b_i
         +\#\big\{b_i\;|\;b_i=0\big\}
         -\sum_{i\not= j}a_ib_j
         \leq
         \sum_{i=1}^m a_i\cdot\sum_{i=1}^mb_i+\sum_{i=1}^m b_i.
       \end{multline*}

     \begin{varthm-roman}[Step 4b]
       If $e<0$, we give several upper bounds for
       $\beta=\sum\limits_{i=1}^m\Big(h^0\big(\Sigma,\ko_\Sigma(\Delta_i)\big)-1\Big)$:
       \begin{displaymath}
         %\sum_{i=1}^m\Big(h^0\big(\Sigma,\ko_\Sigma(\Delta_i)\big)-1\Big)%
         \beta
         \leq
         \as{1.8}
         \left\{
           \begin{array}{ll}
             \tfrac{1}{2}\sum\limits_{i=1}^ma_i \sum\limits_{i=1}^mb_i
             +\tfrac{1}{2}\left(\sum\limits_{i=1}^mb_i\right)^2
             +\tfrac{1}{8}\left(\sum\limits_{i=1}^ma_i\right)^2
             +\tfrac{1}{4}\sum\limits_{i=1}^ma_i
             +\tfrac{1}{2}\sum\limits_{i=1}^mb_i,
             \;\text{if } g=1,\\
             \sum\limits_{i=1}^ma_i \sum\limits_{i=1}^mb_i
             +\sum\limits_{i=1}^ma_i+\sum\limits_{i=1}^mb_i
             -\tfrac{9e}{32}\left(\sum\limits_{i=1}^ma_i\right)^2,
             \;\;\;\text{for } g \text{
               arbitrary.} \\
             \tfrac{1}{4} \sum\limits_{i=1}^ma_i \sum\limits_{i=1}^mb_i
             +\sum\limits_{i=1}^ma_i+\sum\limits_{i=1}^mb_i
             -\tfrac{9e}{32} \left(\sum\limits_{i=1}^ma_i\right)^2
             -\tfrac{1}{2e} \left(\sum\limits_{i=1}^mb_i\right)^2,
             \; g \text{
               arbitrary.} 
           \end{array}
         \right.
         \as{\asf}
       \end{displaymath}
     \end{varthm-roman} 
     If $g$ is arbitrary, the claim follows since 
     a thorough investigation leads to
%     by Corollary
%     \ref{lem:divsors-on-ruled-surfaces} we have
     \begin{displaymath}
       h^0\big(\Sigma,\ko_\Sigma(\Delta_i)\big)\leq
       a_ib_i+a_i+b_i+1-\tfrac{9e}{32}\cdot a_i^2
     \end{displaymath}
     and
     \begin{displaymath}
       h^0\big(\Sigma,\ko_\Sigma(\Delta_i)\big)\leq
       \tfrac{1}{4}\cdot a_ib_i+a_i+b_i+1-\tfrac{9e}{32}\cdot
       a_i^2-\tfrac{1}{2e}\cdot {b_i}^2.
     \end{displaymath}
     If $g=1$, then $e=-1$ and $b=b'+\tfrac{a}{2}$ and we are done since
%     We may once more apply Corollary
%     \ref{lem:divsors-on-ruled-surfaces} and see that in any case
     \begin{multline*}
        h^0\big(\Sigma,\ko_\Sigma(\Delta_i)\big)\leq
        a_ib_i'+b_i'+1+\tfrac{a_i(a_i+1)}{2}+\tfrac{b_i'(b_i'-1)}{2}\\
        =
        \tfrac{1}{2}\cdot a_ib_i+
        \tfrac{1}{2}\cdot {b_i}^2+
        \tfrac{1}{8}\cdot a_i^2+
        \tfrac{1}{4}\cdot a_i+
        \tfrac{1}{2}\cdot b_i+1.
     \end{multline*}

     \begin{varthm-roman}[Step 5]       
       In this last step we gather the information from the previous
       investigations and finish the proof considering a bunch of
       different cases.
     \end{varthm-roman}

     Using Step 3 and Step 4 and taking $\sum_{i=1}^ma_i+b_i\leq m$
     into account, we get the following upper bounds for $\beta'= 
     h^1\big(\Sigma,\kj_{X_0}(D)\big)+
     \sum_{i=1}^m
     \Big(h^0\big(\Sigma,\ko_\Sigma(\Delta_i)\big)-1\Big)$
     \begin{multline*} 
%       h^1\big(\Sigma,\kj_{X_0}(D)\big)+
%       \sum_{i=1}^m
%       \Big(h^0\big(\Sigma,\ko_\Sigma(\Delta_i)\big)-1\Big)\\\leq
       \beta'\leq
       \left\{
         \begin{array}{ll}
           3 \sum\limits_{i=1}^ma_i \sum\limits_{i=1}^mb_i
           +2g  \sum\limits_{i=1}^ma_i,
           &\text{if } e=0,\\
           3 \sum\limits_{i=1}^ma_i \sum\limits_{i=1}^mb_i
           +2g  \sum\limits_{i=1}^ma_i-\tfrac{9e}{32}  \left(\sum\limits_{i=1}^ma_i\right)^2,
           &\text{if } e<0,\\
           \tfrac{9}{4} \sum\limits_{i=1}^ma_i \sum\limits_{i=1}^mb_i
           +2g  \sum\limits_{i=1}^ma_i-\tfrac{9e}{32} 
           \left(\sum\limits_{i=1}^ma_i\right)^2
           -\tfrac{1}{2e} \left(\sum\limits_{i=1}^mb_i\right)^2,
           &\text{if } e<0,\\
           \multicolumn{2}{l}{
           3 \sum\limits_{i=1}^ma_i \sum\limits_{i=1}^mb_i, 
           \hspace*{5.5cm}\text{if } e=0, g=1, \sum\limits_{i=1}^m b_i\not=0,
           }\\
           \multicolumn{2}{l}{
           m\leq\sum\limits_{i=1}^m a_i,
           \hspace*{5.85cm}\text{if } e=0, g=1, \sum\limits_{i=1}^m b_i=0,
           }\\
           \tfrac{5}{2} \sum\limits_{i=1}^ma_i \sum\limits_{i=1}^mb_i
           +\tfrac{1}{2} \left(\sum\limits_{i=1}^mb_i\right)^2
           +\tfrac{1}{8} \left(\sum\limits_{i=1}^ma_i\right)^2
           +\tfrac{5}{4} \sum\limits_{i=1}^ma_i,
           &\text{if } e<0, g=1.
         \end{array}
       \right.
     \end{multline*}
     
     Applying now Step 2b-2d we end up with
     $\tfrac{\beta'\cdot(D-K_\Sigma)^2}{\big(\deg(X_S)\big)^2}\leq\gamma$.
     We thus finally get
     \begin{displaymath}
       \as{1.6}
       \begin{array}{l}
         h^1\big(\Sigma,\kj_{X_0}(D)\big)+
         \sum_{i=1}^m
         \Big(h^0\big(\Sigma,\ko_\Sigma(\Delta_i)\big)-1\Big)
         \;\,=\;\,\beta'\;\,\leq\;\,
         \frac{1}{\gamma\cdot(D-K_\Sigma)^2}\cdot\big(\deg(X_S)\big)^2\\
         \hspace*{1cm}=\;\,
         \frac{1}{\gamma\cdot(D-K_\Sigma)^2}\cdot
         \big(\sum_{z\in\Sigma}\deg(X_{S,z})\big)^2
         \;\,\leq\;\,
         \frac{\# X_S}{\gamma\cdot(D-K_\Sigma)^2}\cdot
         \sum_{z\in\Sigma}\deg(X_{S,z})^2\\
         \hspace*{1cm}\leq\;\,
         \frac{\# X_S}{\gamma\cdot(D-K_\Sigma)^2}\cdot
         \sum_{z\in\Sigma}\deg(X_{0,z})^2\;\,
         <_\expl{(4)}\;\,\# X_S.
       \end{array}
       \as{\asf}
     \end{displaymath}
   \end{proof}

   It remains to show, that the inequality which we derived 
%   in the above cases 
   cannot hold.

   \begin{lemma}\label{lem:irred-E}
     Let $D\in\Div(\Sigma)$, $\ks_1,\ldots,\ks_r$ be pairwise distinct
     topological or analytical singularity types and
     $k_1,\ldots,k_r\in\N\setminus\{0\}$. Suppose that
     $V_{|D|}^{irr,reg}(k_1\ks_1,\ldots,k_r\ks_r)$ is non-empty.

     Then there exists no curve\tom{\footnote{For a subset $U\subseteq V$
       of a topological space $V$ we denote by $\overline{U}$ the
       closure of $U$ in $V$.} }
     $C\in V_{|D|}^{irr}(k_1\ks_1,\ldots,k_r\ks_r)
     \setminus
     \overline{V_{|D|}^{irr,reg}(k_1\ks_1,\ldots,k_r\ks_r)}$
     such that for the zero-dimensional scheme $X_0=X(C)$
     there exist curves 
     $\Delta_1,\ldots,\Delta_m\subset\Sigma$ and
     zero-dimensional locally complete intersections $X_i^0\subseteq
     X_{i-1}$ for $i=1,\ldots,m$, where $X_i=X_{i-1}:\Delta_i$ for
     $i=1,\ldots,m$ %,  such that (a)-(f) in Lemma \ref{lem:irred-A} are fulfilled.
     such that $X_S=\bigcup_{i=1}^mX_i^0$ satisfies
     \begin{equation}
       \label{eq:irred-E:1}
       h^1\big(\Sigma,\kj_{X_0}(D)\big)+
       \sum_{i=1}^m
       \Big(h^0\big(\Sigma,\ko_\Sigma(\Delta_i)\big)-1\Big) < \# X_S.
     \end{equation}
   \end{lemma}

   \begin{proof}
     Throughout the proof we use the notation
     $V^{irr}=V_{|D|}^{irr}(k_1\ks_1,\ldots,k_r\ks_r)$
     and
     $V^{irr,reg}=V_{|D|}^{irr,reg}(k_1\ks_1,\ldots,k_r\ks_r)$.
     
     Suppose there exists a curve $C\in V^{irr}\setminus \overline{V^{irr,reg}}$
     satisfying the assumption of the Lemma, and let $V^*$ be the
     irreducible component of $V^{irr}$
     containing $C$. Moreover, let $C_0\in V^{irr,reg}$.

     We consider in the following the morphism from Subsection \ref{subsec:psi}
     \begin{displaymath}
       \Psi=\Psi_{|D|}(k_1\ks_1,\ldots,k_r\ks_r)
       :V_{|D|}(k_1\ks_1,\ldots,k_r\ks_r)
       \rightarrow
       B(k_1\ks_1,\ldots,k_r\ks_r)=B.
     \end{displaymath}

     \begin{varthm-roman}[Step 1]
       $h^0\big(\Sigma,\kj_{X(C_0)/\Sigma}(D)\big)
       =
       h^0\big(\Sigma,\kj_{X(C)/\Sigma}(D)\big)-
       h^1\big(\Sigma,\kj_{X(C)/\Sigma}(D)\big).$
     \end{varthm-roman}
     By the choice of $C_0$ we have
     \begin{displaymath}
       0=H^1\big(\Sigma,\kj_{X^*(C_0)/\Sigma}(D)\big)\rightarrow
       H^1(\Sigma,\ko_\Sigma(D)\big) \rightarrow
       H^1(\Sigma,\ko_{X^*(C_0)}(D)\big)=0,
     \end{displaymath}
     and thus $D$ is non-special,
     i.~e.~$h^1(\Sigma,\ko_\Sigma(D)\big)=0$.
     But then
     \begin{displaymath}
%       \begin{array}{rcl}
         h^0\big(\Sigma,\kj_{X(C_0)/\Sigma}(D)\big)
%         h^0\big(\Sigma,\ko_\Sigma(D)\big)-\deg\big(X(C_0)\big)\\
%         &=&h^0\big(\Sigma,\ko_\Sigma(D)\big)-\deg\big(X(C)\big)\\
         =
         h^0\big(\Sigma,\kj_{X(C)/\Sigma}(D)\big)-
         h^1\big(\Sigma,\kj_{X(C)/\Sigma}(D)\big).
%       \end{array}
     \end{displaymath}

     \begin{varthm-roman}[Step 2]       
       $h^1\big(\Sigma,\kj_{X(C)}(D)\big)\geq\codim_B\Big(\Psi\big(V^*\big)\Big)$.
     \end{varthm-roman}
     Suppose the contrary, that is
     $\dim\Big(\Psi\big(V^*\big)\Big)<\dim(B)-h^1\big(\Sigma,\kj_{X(C)/\Sigma}(D)\big)$,  
     then by Step 1 and Theorem \ref{thm:v-reg}
     \begin{displaymath}
       \begin{array}{rcl}
         \dim\big(V^*\big)&\leq& \dim\Big(\Psi\big(V^*\big)\Big)+
         \dim\Big(\Psi^{-1}\big(\Psi(C)\big)\Big)\\
         &<&\dim(B)-h^1\big(\Sigma,\kj_{X(C)/\Sigma}(D)\big)+
         h^0\big(\Sigma,\kj_{X(C)/\Sigma}(D)\big)-1\\
         &=&
         \dim(B)+h^0\big(\Sigma,\kj_{X(C_0)/\Sigma}(D)\big)-1=\dim\big(V^{irr,reg}\big).         
       \end{array}
     \end{displaymath}
     However, any irreducible component of $V^{irr}$ has at least the
     expected dimension $\dim\big(V^{irr,reg}\big)$, which gives a
     contradiction. 

     \begin{varthm-roman}[Step 3]
       $\codim_B\Big(\Psi\big(V^*\big)\Big)\geq\#X_S-\sum_{i=1}^m\dim|\Delta_i|_l$.
     \end{varthm-roman}
     The existence of the subschemes $X_i^0\subseteq X(C)\cap
     \Delta_i$ imposes at least $\# X_i^0-\dim|\Delta_i|_l$ conditions
     on $X(C)$ and increases thus the codimension of
     $\Psi\big(V^*\big)$ by this number. 

     \begin{varthm-roman}[Step 4]
       Collecting the results we derive the following contradiction:
     \end{varthm-roman}
     \begin{displaymath}
       \begin{array}{rllll}
         h^1\big(\Sigma,\kj_{X(C)}(D)\big)&\geq_\expl{Step 2}&
         \codim_B\Big(\Psi\big(V^*\big)\Big)\\
         &\geq_\expl{Step 3}&
         \#X_S-\sum_{i=1}^m\dim|\Delta_i|_l
         &>_\expl{\eqref{eq:irred-E:1}}&h^1\big(\Sigma,\kj_{X(C)}(D)\big).
       \end{array}
     \end{displaymath}

%     \begin{varthm-roman}[Step 1]
%       $\dim\big(V^*\big)>\dim\big(V^{irr,reg}\big)$.
%     \end{varthm-roman}
%     Then
%     \begin{displaymath}
%       \begin{array}{rcl}
%         \dim\big(V^*\big)&=&h^0\big(\Sigma,\kj_{X^*(C)/\Sigma}(D)\big)-1\\
%         &=&h^0(\Sigma,\ko_\Sigma(D)\big)-\deg\big(X^*(C)\big)+
%         h^1\big(\Sigma,\kj_{X^*(C)/\Sigma}(D)\big)-1\\
%         &>&h^0(\Sigma,\ko_\Sigma(D)\big)-\deg\big(X^*(C_0)\big)-1
%         =\dim\big(V^{irr,reg}\big).
%       \end{array}
%     \end{displaymath}
   \end{proof}

   \bibliographystyle{amsalpha-tom}
   \bibliography{bibliographie}

\end{document}